\documentclass[10pt,journal]{IEEEtran}

\usepackage{amsfonts}
\usepackage{extarrows}
\usepackage{enumerate}
\usepackage{amsmath}
\usepackage{amssymb}
\usepackage{flushend}
\usepackage{graphicx}
\usepackage{color}
\usepackage[square,comma,sort&compress,numbers]{natbib}

\allowdisplaybreaks[4]

\begin{document}

\title{Reinforcement Learning-based Disturbance Rejection Control for Uncertain Nonlinear Systems}


\author{Maopeng Ran, Juncheng Li, and Lihua Xie, \emph{Fellow, IEEE}
 \thanks{This work is supported by the A$^*$STAR Industrial Internet of Things Research Program, under the RIE2020 IAF-PP Grant A1788a0023.}
\thanks{M. Ran, J. Li, and L. Xie (corresponding author) are with the School of Electrical and Electronic Engineering, Nanyang Technological University, Singapore 639798 (email: mpran@ntu.edu.sg; juncheng001@e.ntu.edu.sg; elhxie@ntu.edu.sg).}}

\maketitle

\begin{abstract}

This paper investigates the reinforcement learning (RL) based disturbance rejection control for uncertain nonlinear systems having non-simple nominal models. An extended state observer (ESO) is first designed to estimate the system state and the total uncertainty, which represents the perturbation to the nominal  system dynamics. Based on the output of the observer, the control compensates for the total uncertainty in real time, and simultaneously, online approximates the optimal policy for the compensated system using a simulation of experience based RL  technique. Rigorous theoretical analysis is given to show the practical convergence of the system state to the origin and the developed policy to the ideal optimal policy. It is worth mentioning that, the widely-used restrictive persistence of excitation (PE) condition is not required in the established framework. Simulation results are presented to illustrate the effectiveness of the proposed method.

\end{abstract}


\IEEEpeerreviewmaketitle

\section{Introduction}


Reinforcement learning (RL), inspired by learning behavior in nature, is a goal-oriented learning strategy wherein the agent learns the policy to optimize a pre-defined reward by interacting with the environment \cite{Nef-2019,Lewis-2009}. From a perspective of control theory, RL is strongly connected with traditional optimal control and adaptive control algorithms \cite{Lewis-2018}. The goal of an RL-based controller is to learn the optimal policy and value function by finding the solution to the Hamilton-Jacobi-Bellman (HJB) equation online. For its uniqueness in data-driven concepts, effectiveness in reaching optimal behavior, and adaptiveness to uncertain environment, RL has undergone rapid progress in control community \cite{Liu-2017,Lewis-2018}.

Early-efforts on RL-based control were mainly devoted to linear systems. For example, a policy iteration technique to solve the continuous time LQR problem without the knowledge of the  state matrix was proposed in \cite{Lewis-2009b}. The results in \cite{Lewis-2009b} were generalized to linear systems with completely unknown dynamics in  \cite{Jiang-2012}. In the last decade, the extension of RL-based control to nonlinear systems has gained considerable attention as well. Model-free optimal control problems for affine and nonaffine nonlinear systems were addressed based on RL in \cite{Luo-2014} and \cite{Jiang-2014}, respectively. In \cite{Jiang-2015b}, RL-based optimal control for polynomial systems was considered. A new policy iteration scheme was developed in \cite{Jiang-2015b} which attempts to solve a convex optimization problem rather than a partial differential equation as usual. In \cite{Zarg-2015} and \cite{Jiang-2017}, the problem of RL-based optimal tracking control for nonlinear systems in strict-feedback form was investigated. Ref. \cite{Zarg-2015} handled this problem by transforming it into an equivalent optimal regulation problem through a feedforward adaptive input, while \cite{Jiang-2017} addressed this problem by using a novel policy iteration technique to solve semidefinite HJB equations.

In practice, almost all control systems are inherently nonlinear and subject to uncertainties and disturbances \cite{Khalil-2002,Isidori-1989}. Thus it is important to investigate the robustness of RL-based control policies for nonlinear systems. The traditional robust control methods, such as the nonlinear small-gain theorem \cite{Jiang-2014b,Gao-2019}, sliding mode control \cite{Fan-2016}, and $H_{\infty}$ control \cite{Xue-2020,Jiang-2017b} were integrated with the RL designs. In \cite{Wang-2018,Wang-2020}, novel adaptive critic strategies were developed with robustness guarantees for disturbed nonlinear systems.  In \cite{Zhang-2018,Yang-2019b}, event-triggered mechanisms were introduced into the design of robust RL controllers to mitigate unnecessary communications. Note that to guarantee parameter convergence, the RL-based control laws generally require persistence of excitation (PE) of the system state.  In \cite{Jiang-2014b,Fan-2016,Wang-2020,Wang-2018}, to fulfill the PE condition, probing signals are added into the control inputs, which will inevitably decrease the transient performance. What is more,  the approaches in \cite{Jiang-2014b,Gao-2019,Fan-2016,Xue-2020,Jiang-2017b,Wang-2018,Wang-2020,Zhang-2018,Yang-2019b} were developed only for relatively simple uncertainties or disturbances (e.g., independent from the control signal).




On the other hand, the active disturbance rejection control (ADRC), which was proposed by Han \cite{Han-2009}, is an efficient methodology to deal with uncertainties and disturbances. The basic idea of ADRC is to first estimate the total uncertainty online via a device called extended state observer (ESO), and then compensate for it in the control loop in real time. Accordingly, an ADRC law contains two components: a disturbance rejection term compensates for the estimated uncertainty, and a nominal controller guarantees the performance of the compensated system. The existing ADRC literature mainly focus on its engineering applications \cite{Sari-2020} and theoretical verifications for different classes of uncertain systems, such as MIMO systems \cite{Guo-2013}, nonaffine-in-control systems \cite{Ran-2017a}, time-delay systems \cite{Ran-2020}, and so on. On the contrary, there are very few works that concern the design of the nominal controller, since the compensated system is generally simplified into a chain of integrators (see, e.g., \cite{Guo-2013,Ran-2017a,Ran-2020,Gao-2003}).

Motivated by the observations stated above, in this paper, we consider the RL-based disturbance rejection control for uncertain nonlinear systems. The systems are assumed to have non-simple nominal models (i.e., the nominal models maybe also complex and nonlinear) and are subject to multiple uncertainties. In order to relax the requirement of the PE condition and avoid using the probing signal, concurrent learning (CL), also known as experience replay, is employed in this paper. The original idea of concurrent learning is to use the recorded data and current data simultaneously for parameter adaptation in the framework of model reference adaptive control (MRAC) \cite{Chowdhary-2010}. In recent years, this idea was extended to develop RL-based control policies to remove the PE condition \cite{Zhang-2018,Xue-2020,Yang-2019b,Yang-2020a,Vam-2015}. Note that in RL, estimates for the parameters of the value function are updated using the Bellman error (BE) as a performance metric. In this paper, we evaluate the BE along the system trajectory,  and simultaneously, extrapolate the BE to any unexplored data points by using the nominal model. This concurrent learning technique is interpreted as simulation of experience \cite{Kama-2016a,Kama-2016b}.
The main contributions of this paper are twofold:
\begin{enumerate}
  \item An ESO-based reinforcement learning and disturbance rejection framework for uncertain nonlinear systems having non-simple nominal models is established. The ESO provides estimates of the system state and total uncertainty to the RL-algorithm and disturbance rejection term, respectively.
      The disturbance rejection term compensates for the total uncertainty in real time, while the    output feedback simulation of experience based RL explores the optimal policy for the compensated system simultaneously. To the best of the authors' knowledge, this work is the first attempt that bridges the gap between RL and the philosophy of ADRC.
  \item A rigorous theoretical analysis that establishes online simultaneous disturbance rejection and optimal policy approximation without the PE condition  is conducted. Due to the strong coupling between the ESO estimation error and the RL approximation error, the theoretical problem faced in this paper is much more challenging than the previous ADRC \cite{Guo-2013,Ran-2017a,Ran-2020,Gao-2003} and RL results \cite{Kama-2016a,Kama-2016b}. Specifically, in \cite{Guo-2013,Ran-2017a,Ran-2020,Gao-2003} the nominal controller is designed to stabilize a chain of integrators, while here it is a learning controller which is able to online approximate the optimal policy for the non-simple compensated system. Compared with \cite{Kama-2016a,Kama-2016b}, we provide a more practical RL solution for nonlinear systems, since the developed controller is capable of handling multiple uncertainties and it is output feedback based.
\end{enumerate}

The rest of this paper is organized as follows. Section II presents some preliminaries and formulates the problem. Section III states the design and analysis of the developed ESO-based reinforcement learning and disturbance rejection framework. In Section IV, some simulation results are given to demonstrate the effectiveness of the proposed method. Finally, Section V concludes the paper.

\emph{Notations:} Throughout this paper, we use $\mathcal{C}$ to represent the set of all continuously differentiable functions. For any continuously differentiable function $f: \mathbb{R}^{n}\times \mathbb{R}^m\rightarrow \mathbb{R}^l$, $f_x:\mathbb{R}^{n}\times \mathbb{R}^m\rightarrow \mathbb{R}^{l\times n} $ represents its gradient with respect to the first vector argument, i.e., $f_x(\nu_1, \nu_2)=\partial f(\nu_1,\nu_2)/\partial \nu_1$.  $\lambda_{\max}(P)$ and $\lambda_{\min}(P)$ denote the maximum and minimum eigenvalues of matrix $P$, respectively. Big $O$-notation in terms of $\nu$ is denoted as $O(\nu)$ and it is assumed that this holds for $\nu$ positive and sufficiently small. $\textrm{sat}:\mathbb{R}\rightarrow \mathbb{R}$ denotes the standard unity saturation function defined by $\textrm{sat}(\nu)=\textrm{sign}(\nu)\cdot\min\{1, |\nu|\}$. $\textbf{1}_A(\nu)$ is the indictor function defined by $\textbf{1}_A(\nu)=\left\{\begin{matrix} 1 ~~ \textrm{if} ~ \nu\in A, \\ 0 ~~ \textrm{if} ~ \nu \notin A. \end{matrix}\right.$ Throughout this paper, for the sake of brevity, the time variable $t$ of a signal will be omitted except when the dependence of the signal on $t$ is crucial for clear presentation.

\section{Problem Formulation}

A single-input, single-output system with relative degree $n$, under a suitable diffeomorphism, can be written in the following normal form \cite{Isidori-1989}:
\begin{equation}\label{eq1}
 \left\{
  \begin{aligned}
          \dot{z}=& f_z(x, z, \omega), \\
          \dot{x}=& Ax+B[f(x, z, \omega)+g(x,z,\omega)u], \\
                y=& Cx,
        \end{aligned} \right.
\end{equation}
where $x=[x_1, \ldots, x_n]^{\rm{T}}\in\mathbb{R}^n$ and $z\in\mathbb{R}^p$ are the states, $\omega\in\mathbb{R}$ is the external disturbance, $u\in\mathbb{R}$ is the control input, $y\in\mathbb{R}$ is the measured output, $f_z: \mathbb{R}^n\times \mathbb{R}^p \times \mathbb{R}\rightarrow \mathbb{R}^p$ and  $f, g: \mathbb{R}^n\times \mathbb{R}^p \times \mathbb{R}\rightarrow \mathbb{R}$ are continuously differentiable functions,  and matrices $A\in\mathbb{R}^{n\times n}$, $B\in\mathbb{R}^{n\times 1}$, and $C\in\mathbb{R}^{1\times n}$ are given by
\begin{equation*}
   A=
\left[
  \begin{array}{ccccc}
    0 & 1 &  \cdots & 0 \\
    \vdots &  \vdots & \ddots & \vdots \\
    0 & 0 &  \cdots & 1 \\
    0 & 0 &    \cdots  &0 \\
  \end{array}
\right], ~B=\left[
             \begin{array}{c}
               0 \\
               \vdots \\
               0 \\
               1 \\
             \end{array}
           \right],~C=\left[
             \begin{array}{c}
               1 \\
               0 \\
               \vdots \\
               0 \\
             \end{array}
           \right]^{\rm{T}}.
 \end{equation*} For system (\ref{eq1}), the following assumptions are made:

\emph{Assumption A1:} The external disturbance $\omega$ and its time derivative $\dot{\omega}$  are bounded.

\emph{Assumption A2:} The zero-dynamics $\dot{z}=f_z(x, z, \omega)$ with input $(x, \omega)$ is bounded-input-bounded-state stable.

In this paper, the functions $f(\cdot)$ and $g(\cdot)$ are assumed to be uncertain and partially known, that is,
\begin{align}\label{eq2}
  f(x,z,\omega)=& f_0(x)+\Delta f(x,z,\omega), \\ ~g(x,z,\omega)=& g_0(x)+\Delta g(x,z,\omega),
\end{align}
where $f_0(\cdot), g_0(\cdot)\in \mathcal{C}(\mathbb{R}^n, \mathbb{R})$ are known and globally bounded, $g_0(\cdot)\neq 0$, and $\Delta f(\cdot), \Delta g(\cdot)\in\mathcal{C}(\mathbb{R}^n\times\mathbb{R}^p\times\mathbb{R}, \mathbb{R})$ are unknown.

\emph{Remark 1:} Similar to many previous ESO works (see, e.g., \cite{Khalil-2010,Khalil-2008}), we require the global boundedness of $f_0(\cdot)$ and $g_0(\cdot)$. It should be pointed out that this requirement does not exclude linear functions or any unbounded function
because in this paper we achieve semi-global results, and the global boundedness can be always satisfied by smoothly saturating the function outside a compact set of interest. \IEEEQED

\emph{Remark 2:} In system (\ref{eq1}), the uncertainties are located in the same channel as the control signal. For a general nonlinear system with smooth unmatched uncertainties, one can conduct a well-defined state transformation to transform the system into the form of system (\ref{eq1}) \cite{Zhao-2017}. \IEEEQED

Note that system (\ref{eq1}) is subjected to multiple uncertainties, including the uncertainty in  zero dynamics $\dot{z}=f_z(x, z, \omega)$, drift dynamics modeling error $\Delta f(\cdot)$,  mismatch of control $\Delta g(\cdot)u$, and external disturbance $\omega$. Inherited from the basic philosophy of ADRC \cite{Han-2009}, the total uncertainty is regarded as an extended state of the system, which is denoted by
\begin{equation}\label{eq8}
  x_{n+1}\triangleq \Delta f(x,z,\omega)+\Delta g(x,z,\omega)u.
\end{equation}

With assumption A2, we concentrate on the control design for the $x$-subsystem.
Had $x$ and $x_{n+1}$ been available for feedback, we could have the following controller:
\begin{equation}\label{eq9}
  u=u_0^*(x)-\frac{x_{n+1}}{g_0(x)}.
\end{equation}
Here the second term $-\frac{x_{n+1}}{g_0(x)}$ is to compensate the total uncertainty, and $u_0^*(x)$ is the optimal policy for the compensated system
\begin{equation}\label{eq23}
  \dot{x}=Ax+B[f_0(x)+g_0(x)u_0],
\end{equation}
with the cost functional
\begin{equation}\label{eq3}
  J(x,u_0)=\int_{0}^{\infty}r(x(\tau),u_0(\tau))\textrm{d}\tau.
\end{equation}
The function $r: \mathbb{R}^n\times \mathbb{R}\rightarrow \mathbb{R}^+$ denotes the instantaneous cost given by
\begin{equation}\label{eq24}
  r(x,u_0)=Q(x)+u_0^{\textrm{T}}Ru_0,
\end{equation}
where $Q:\mathbb{R}^n\rightarrow \mathbb{R}^+$ is positive definite, and $R>0$.

Note that in the control (\ref{eq9}), the optimal control  $u_0^*(x)$ is generally analytically infeasible, and the state $x$ and total uncertainty $x_{n+1}$ are both unavailable. The objective of this paper is then to develop an ESO-based reinforcement learning and disturbance rejection scheme, which is able to online simultaneously compensate the uncertainty  and approximate the optimal policy.

\section{ESO-based Reinforcement Learning and Disturbance Rejection}


\subsection{Control Design}

First of all, an ESO is designed for system (\ref{eq1}) to provide the estimates of the state and total uncertainty:
\begin{equation}\label{eso}
 \left\{
  \begin{aligned}
& \dot{\widehat{x}}_1=\widehat{x}_2+\frac{l_1}{\varepsilon}(x_1-\widehat{x}_1),   \\
& ~\vdots \\
&\dot{\widehat{x}}_{n-1}=\widehat{x}_n+\frac{l_{n-1}}{\varepsilon^{n-1}}(x_1-\widehat{x}_1),  \\
&\dot{\widehat{x}}_{n}=\widehat{x}_{n+1}+\frac{l_n}{\varepsilon^n}(x_1-\widehat{x}_1)+f_0(\widehat{x})+g_0(\widehat{x})u, \\
&\dot{\widehat{x}}_{n+1}=\frac{l_{n+1}}{\varepsilon^{n+1}}(x_1-\widehat{x}_1),
  \end{aligned}
\right.
\end{equation}
where $\widehat{x}=[\widehat{x}_1,\ldots,\widehat{x}_{n}]^{\rm{T}}\in\mathbb{R}^{n}$, $\varepsilon<1$ is a small positive constant, $L=[l_1, l_2,\ldots,l_{n+1}]^{\rm{T}}\in\mathbb{R}^{n+1}$ is the observer gain selected such that the following matrix is Hurwitz:
\begin{equation*}\label{eq14}
  E=\left[
      \begin{array}{ccccc}
        -l_1 & 1 & 0 & \cdots & 0 \\
        -l_2 & 0 & 1 & \cdots & 0 \\
        \vdots &\vdots &\vdots & \ddots & \vdots \\
        -l_n & 0 & 0 & \cdots & 1 \\
        -l_{n+1} & 0 & 0 & \cdots & 0 \\
      \end{array}
    \right]\in\mathbb{R}^{(n+1)\times(n+1)}.
\end{equation*}
Since $\varepsilon$ is a small positive constant, the observer (\ref{eso}) exhibits peaking phenomenon during the initial fast transient \citep{Guo-2013,Khalil-2008}. To protect the peaking of the observer from propagating into other variables, we employ the well-known saturation technique \cite{Khalil-2008} to saturate the output of the observer as
\begin{equation}\label{eq28}
\overline{x}_i=M_is\left(\frac{\widehat{x}_i}{M_i}\right), ~1\leq i\leq n+1,
\end{equation}
where $M_i$, $1\leq i\leq n+1$, are saturation bounds to be selected such that the saturation will not be invoked in the steady-period of the observer, and $s(\cdot)$ is an odd saturation-like function defined by
\begin{equation*}
s(\nu)= \left\{
  \begin{aligned}
       &   \nu, && 0\leq \nu \leq 1 \\
       &   \nu+\frac{\nu-1}{\varepsilon}-\frac{\nu^2-1}{2\varepsilon}, && 1\leq \nu \leq 1+\varepsilon \\
      &          1+\frac{\varepsilon}{2}, &&  \nu>1+\varepsilon
        \end{aligned} \right.
\end{equation*}
Note that this function is nondecreasing, continuously differentiable with locally Lipschitz derivative. What is more, $0\leq s'(\nu)\leq 1$ and $\left|s(\nu)-\textrm{sat}(\nu)\right|\leq \frac{\varepsilon}{2}$, $\forall \nu\in\mathbb{R}$, where $s'(\nu)\triangleq \frac{\textrm{d}s(\nu)}{\textrm{d}\nu}$. For subsequent use, denote $\overline{x}=[\overline{x}_1,\ldots, \overline{x}_n]^{\rm{T}}$, and $\dot{\overline{x}}_i=s'(\widehat{x}_i/M_i)\dot{\widehat{x}}_i$, $1\leq i\leq n+1$.

According to the optimal control theory \cite{Kirc-2012}, the optimal control problem of system (\ref{eq23}) with cost functional (\ref{eq3}) can be converted to solve the following HJB equation:
\begin{align}\label{eq5}
  {V^*_x}(x)\left[Ax+B(f_0(x)+g_0(x)u_0^*(x))\right] \qquad  \nonumber \\
  +Q(x)+{u_0^*}^{\textrm{T}}(x)Ru_0^*(x)=0,
\end{align}
where $V^*\in\mathcal{C}(\mathbb{R}^n,\mathbb{R}_{\geq 0})$, $V^*(0)=0$, is the optimal value function. The optimal control policy can be determined from the optimal value function as
\begin{equation}\label{eq6}
  u_0^*(x)=-\frac{1}{2}R^{-1}g_0^{\rm{T}}(x)B^{\rm{T}}V^{*}_x{^{\rm{T}}}(x).
\end{equation}
Generally, the analytical solution of the HJB equation is not feasible. However,  the optimal value function $V^*(x)$ and the optimal control policy $u_0^*(x)$ can be approximated by an actor-critic neural network (NN) based approach  \cite{Jiang-2014,Luo-2014}. For any given compact set $\mathcal{X}\subset \mathbb{R}^n$ and positive constant $\overline{\kappa}$, the optimal value function $V^*(x)$ can be represented by
\begin{equation}\label{eq50}
  V^*(x)=\Theta^{\textrm{T}}\phi(x)+\kappa(x),
\end{equation}
where $\phi: \mathcal{X}\rightarrow \mathbb{R}^{l}$ is a continuously differentiable activation function, $\Theta\in\mathbb{R}^l$ is the ideal weight vector, $l\in\mathbb{N}$ is the number of neurons, and $\kappa:\mathbb{R}^n\rightarrow \mathbb{R}$ is the approximation error function which satisfies $|\kappa(x)|\leq \overline{\kappa}$ and $|\kappa_x(x)|\leq \overline{\kappa}$, $\forall x\in\mathcal{X}$. It follows from the NN representation of the value function that the  optimal control policy can be obtained as
\begin{equation}\label{eq51}
  u_0^*(x)=-\frac{1}{2}R^{-1}g_0^{\rm{T}}(x)B^{\rm{T}}\left(\phi_x^{\textrm{T}}(x)\Theta+\kappa^{\rm{T}}_x(x)\right).
\end{equation}
Therefore, bearing in mind that the ESO (\ref{eso}) provides an estimate of the state $x$, the NN-based approximations of $V^*(x)$ and $u_0^*(x)$ are given by
\begin{align}
\label{eq25}   \widehat{V}\left(\overline{x},\widehat{\Theta}_v\right)= & \widehat{\Theta}^{\textrm{T}}_v\phi(\overline{x}), \\
\label{eq26}  \widehat{u}_0\left(\overline{x},\widehat{\Theta}_c\right)= & -\frac{1}{2}R^{-1}g_0^{\rm{T}}(\overline{x})B^{\rm{T}}\phi^{\rm{T}}_x(\overline{x})\widehat{\Theta}_c,
\end{align}
where $\widehat{\Theta}_v, \widehat{\Theta}_c\in\mathbb{R}^l$ are estimates of $\Theta$, and weights for the critic and actor NNs, respectively.  By (\ref{eq25}) and (\ref{eq26}), the approximated instantaneous BE $\delta_t:\mathbb{R}^{n}\times \mathbb{R}^l\times \mathbb{R}^l\rightarrow \mathbb{R}$ can be computed as
\begin{align}\label{eq7}
\delta_t \triangleq &
  \delta\left(\overline{x},\widehat{\Theta}_v,\widehat{\Theta}_c\right) \nonumber \\
  \triangleq  &\widehat{V}_x(\overline{x},\widehat{\Theta}_v)\left[A\overline{x}+B\left(f_0(\overline{x})+g_0(\overline{x})\widehat{u}_0(\overline{x},\widehat{\Theta}_c)\right)\right]\nonumber \\
  & +Q(\overline{x})+\widehat{u}_0^{\textrm{T}}(\overline{x},\widehat{\Theta}_c)R\widehat{u}_0(\overline{x},\widehat{\Theta}_c).
\end{align}

In this paper, motivated by \cite{Kama-2016a,Kama-2016b}, we implement an ESO based learning strategy via simulation of experience, i.e., use the estimated state $\widehat{x}$ and the knowledge of $f_0(\cdot)$ and $g_0(\cdot)$ to extrapolate the BE to a predefined set of points $\{x^i\in\mathbb{R}^n|i=1,\ldots, N\}$. The approximated BE extrapolated to the point $x^i$ is given by
\begin{equation}\label{eq49}
  \delta_i\triangleq\delta\left(x^i,\widehat{\Theta}_v,\widehat{\Theta}_c\right).
\end{equation}
Then the actor-critic NN uses the BEs $\delta_t$ and $\delta_i$ to update the estimates $\widehat{\Theta}_v$ and $\widehat{\Theta}_c$. A least-square update law for the critic NN is given by
\begin{align}
\label{eq15}  \dot{\widehat{\Theta}}_v=& -\lambda_{v1}\Gamma\frac{\mu}{\rho}\delta_t-\frac{\lambda_{v2}}{N}\Gamma \sum_{i=1}^{N}\frac{\mu_i}{\rho_i}\delta_i, \\
\label{eq16}  \dot{\Gamma}= & \left(\beta\Gamma-\lambda_{v1}\frac{\Gamma\mu\mu^{\rm{T}}\Gamma}{\rho^2}\right)\textbf{1}_{\{\|\Gamma\|\leq \varsigma_1\}}, ~\|\Gamma(0)\|\leq \varsigma_1,
\end{align}
where $\Gamma:\mathbb{R}_{\geq 0}\rightarrow \mathbb{R}^{l\times l}$ is a time-varying least-square gain matrix, $\varsigma_1>0$ is a saturation constant,  $\lambda_{v1}, \lambda_{v2}>0$ are constant adaption gains, $\beta>0$ is a constant forgetting factor, and
\begin{align*}
  \mu = & \phi_x(\overline{x})\left[A\overline{x}+B\left(f_0(\overline{x})+g_0(\overline{x})\widehat{u}_0(\overline{x},\widehat{\Theta}_c)\right)\right], \\
  \mu_i = & \phi_x(x^i)\left[Ax^i+B\left(f_0(x^i)+g_0(x^i)\widehat{u}_0(x^i,\widehat{\Theta}_c)\right)\right],\\
  \rho = & 1+\gamma\mu^{\rm{T}}\Gamma\mu, \\
  \rho_i = & 1+\gamma\mu_i^{\rm{T}}\Gamma\mu_i,
\end{align*}
with $\gamma>0$ a constant normalization gain. The update law (\ref{eq16}) indicates that the time-varying gain matrix $\Gamma$ is bounded in the sense that \cite{Kama-2016b}
\begin{equation}\label{eq54}
  \varsigma_0I \leq \Gamma(t)\leq \varsigma_1 I, ~t\geq 0,
\end{equation}
where $\varsigma_0\in\mathbb{R}^+$ and $I$ is the identity matrix. According to the subsequent stability analysis, the update law for the actor NN is given by
\begin{align}\label{eq17}
 \dot{\widehat{\Theta}}_c= & -\lambda_{c1}\left(\widehat{\Theta}_c-\widehat{\Theta}_v\right)-\lambda_{c2}\widehat{\Theta}_c\nonumber \\
& +\frac{\lambda_{v1}G_{t}^{\rm{T}}\widehat{\Theta}_c\mu^{\textrm{T}}}{4\rho}\widehat{\Theta}_v + \sum_{i=1}^{N}\frac{\lambda_{v2}G^{\rm{T}}_{i}\widehat{\Theta}_c\mu^{\rm{T}}_i}{4N\rho_i}\widehat{\Theta}_v,
\end{align}
where $\lambda_{c1},\lambda_{c2}>0$ are constant adaption gains, and
\begin{align*}
  G_{t}\triangleq & \phi_x(\overline{x})Bg_0(\overline{x})R^{-1}g_0^{\rm{T}}(\overline{x})B^{\rm{T}}\phi^{\rm{T}}_x(\overline{x}),\\
  G_{i}\triangleq & \phi_x(x^i)Bg_0(x^i)R^{-1}g^{\rm{T}}_0(x^i)B^{\rm{T}}\phi_x^{\rm{T}}(x^i).
\end{align*}

Finally, based on the output of the ESO (\ref{eso}), and (\ref{eq26}), the control that simultaneously compensates the total uncertainty and approximates the optimal policy is given by
\begin{equation}\label{eq27}
  u=\widehat{u}_0\left(\overline{x},\widehat{\Theta}_c\right)-\frac{\overline{x}_{n+1}}{g_0(\overline{x})}.
\end{equation}

The block diagram of the developed ESO-based reinforcement learning and disturbance rejection is depicted in Fig. 1.

\begin{figure}
  \centering
  \includegraphics[width=0.5\textwidth,bb=10 10 490 220, clip]{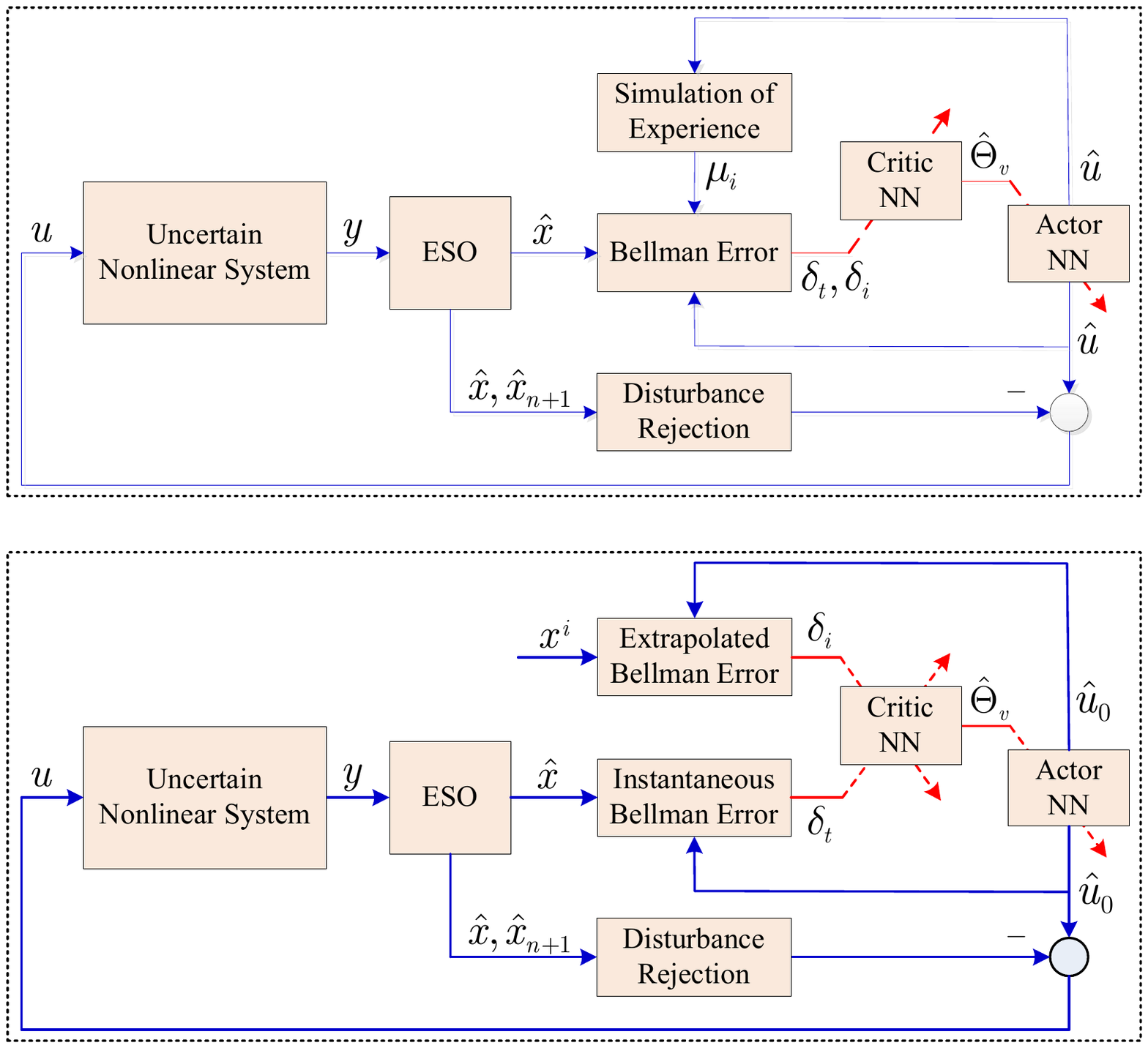}
  \caption{ESO-based reinforcement learning and disturbance rejection architecture.}
\end{figure}


\subsection{Stability Analysis}

To establish the stability of the closed-loop system formed by (\ref{eq1}), (\ref{eso}), (\ref{eq15}), (\ref{eq16}), (\ref{eq17}), and (\ref{eq27}),  we need two additional assumptions. By Assumption A1, let $\omega(t)\in \mathcal{W}$, $\forall t\in[0,\infty)$, for some compact set $\mathcal{W} \subseteq \mathbb{R}$. Assume that the initial condition of system (\ref{eq1}), $x(0)$, belongs to a compact set $\mathcal{X}_0$ which is slightly smaller than $\mathcal{X}$ (i.e., $\mathcal{X}_0\subseteq \mathcal{X}$ and their boundaries are disjoint). By Assumption A2, there exists a constant $c_z$ such that the compact set $\mathcal{Z}=\{z\in\mathbb{R}^p; \|z\|\leq c_z\}$ is a positive invariant set of $\dot{z}=f_z(x,z,\omega)$, for all $x\in\mathcal{X}$ and $\omega \in \mathcal{W}$. Denote the compact set $\Xi$ by $\Xi=\mathcal{X}\times \mathcal{Z}\times \mathcal{W}$.

\emph{Assumption A3:} The known control coefficient $g_0(x)$ satisfies
\begin{equation}\label{eq29}
  k_g\triangleq \max\limits_{(x,z,\omega)\in \Xi, \overline{x}\in\mathcal{X}}\left|\frac{g(x,z,\omega)-g_0(x)}{g_0(\overline{x})}\right|<1.
\end{equation}

\emph{Assumption A4:} There exists a finite set of fixed points $\{x^i\in\mathbb{R}^n|i=1,\ldots, N\}$ such that,
\begin{equation}\label{eq30}
  c \triangleq  \frac{1}{N}\inf\limits_{t\geq 0}\left(\lambda_{\min}\left\{\sum_{i=1}^N\frac{\mu_i\mu_i^{\rm{T}}}{\rho_i}\right\}\right)>0.
\end{equation}

\emph{Remark 3:} Assumption A3 is a standard assumption for the ESO-based controller \cite{Khalil-2008,Guo-2013}. The rational behind this assumption is that the total uncertainty $x_{n+1}$ to be compensated in the control law contains the mismatch of control (i.e., $\Delta g(x,z,\omega)u$). To guarantee the ``compensation capability'' of the ESO-based controller, the nominal control coefficient $g_0(x)$ should not be too far away from the actual control coefficient $g(x,z,\omega)$. Condition (\ref{eq29}) is not difficult to be satisfied in practice \cite{Khalil-2008}. Specifically, for any uncertain control gain $g(x,z,\omega)$ satisfying $\underline{g}<g(x,z,\omega)|_{(x,z,\omega)\in \Xi}<\overline{g}$ (or $-\overline{g}<g(x,z,\omega)|_{(x,z,\omega)\in \Xi}<-\underline{g})$, where $\overline{g}$ and $\underline{g}$ are positive constants. By selecting $g_0(\overline{x})|_{\overline{x}\in\mathcal{X}}>\overline{g}/2$ (or $g_0(\overline{x})|_{\overline{x}\in\mathcal{X}}<-\overline{g}/2)$, the condition (\ref{eq29}) is satisfied. \IEEEQED

\emph{Remark 4:} Condition (\ref{eq30}) depends on the estimate $\widehat{\Theta}_c$, and hence it cannot be guaranteed a prior. However, different from the well-known PE condition \cite{Tao-2003}, the condition (\ref{eq30}) can be monitored online \cite{Kama-2016a}. Moreover, this condition can be heuristically met by selecting more points than the number of the neurons, i.e., choosing $N\gg l$ (see the numerical example in Section IV). \IEEEQED

Now, we are in a position to state our main result.

\emph{Theorem 1:} Consider the closed-loop system formed by plant (\ref{eq1}), ESO (\ref{eso}), control (\ref{eq27}), and RL-based update laws (\ref{eq15}), (\ref{eq16}), and (\ref{eq17}). Suppose Assumptions A1 to A4 are satisfied, and the initial condition of the plant $x(0)$ is an interior point of $\mathcal{X}_0$. Then there exists $\varepsilon^*>0$ such that for any $\varepsilon\in(0,\varepsilon^*)$:
\begin{itemize}
  \item for any $T>0$, $\lim_{\varepsilon\rightarrow 0}|x_i-\widehat{x}_i|\rightarrow 0$, $1\leq i\leq n+1$, uniformly in $t\in[T,\infty)$;
  \item the state $x$ and the weight estimation errors $\widetilde{\Theta}_v\triangleq \widehat{\Theta}_v-\Theta$ and $\widetilde{\Theta}_c\triangleq \widehat{\Theta}_c-\Theta$ are uniformly ultimately bounded.
\end{itemize}

\emph{Proof:} First of all,  we  give some notations and definitions for subsequent use.  Let $Z(t)=[x^{\rm{T}}(t), ~\widetilde{\Theta}_v^{\rm{T}}(t),$ $\widetilde{\Theta}_c^{\rm{T}}(t)]^{\rm{T}}$, $V_v=\frac{1}{2}\widetilde{\Theta}_v^{\rm{T}}\Gamma^{-1}\widetilde{\Theta}_v$, and $V_c=\frac{1}{2}\widetilde{\Theta}_c^{\rm{T}}\widetilde{\Theta}_c$. Denote $\tau_v=\frac{1}{2}\varsigma_1\|\widetilde{\Theta}_v(0)\|+1$, $\tau_c=$ $\max\{V_c(\widetilde{\Theta}_c(0)),V_c(\widetilde{\Theta}_c)_{\|\Theta_c\|\geq \iota_7/|\iota_6|}\}+1$, where $\iota_6$ and $\iota_7$ are given by (\ref{eq46}).  Define several compact sets:
\begin{align*}
\Omega_{v}^0=&\{\widetilde{\Theta}_v\in\mathbb{R}^l; V_v(\widetilde{\Theta}_v)\leq \tau_v\}, \\
\Omega_{v}^1=&\{\widetilde{\Theta}_v\in\mathbb{R}^l; V_v(\widetilde{\Theta}_v)\leq \tau_v+1\}, \\
\Omega_{c}^0=&\{\widetilde{\Theta}_c\in\mathbb{R}^l; V_c(\widetilde{\Theta}_c)\leq \tau_c\}, \\
\Omega_{c}^1=&\{\widetilde{\Theta}_c\in\mathbb{R}^l; V_c(\widetilde{\Theta}_c)\leq \tau_c+1\}.
\end{align*}
Based on the above definitions, denote $\Omega^0=\mathcal{X}_0\times\Omega^0_v\times \Omega^0_c$ and $\Omega^1=\mathcal{X}\times\Omega^1_v\times \Omega^1_c$. It can be clearly observed that $Z(0)$ is an internal point of $\Omega^0$ and $\Omega^0\subseteq \Omega^1-\partial \Omega^1$.

 The following proof will be divided into two steps. In the first step, we show that for sufficiently small $\varepsilon$ and appropriately selected learning gains,  $Z(t)\in\Omega^1$, $\forall t\in [0, \infty)$. Then, in the second step, the convergence of the ESO and the uniformly ultimately boundedness of $x$, $\widetilde{\Theta}_v$, and $\widetilde{\Theta}_c$  are proved.

\emph{Step 1): There exists $\varepsilon^{\dag}>0$ such that for any $\varepsilon\in(0, \varepsilon^{\dag})$, $Z(t)\in\Omega^1$, $\forall t\in [0, \infty)$.}

Since $Z(0)$ is an interior point of $\Omega^0$, $\Omega^0\subseteq \Omega^1-\partial \Omega^1$, and the output of the ESO (\ref{eso}) is saturated, there exists an $\varepsilon$-independent $t_0>0$ such that $Z(t)\in\Omega^1$, $\forall t\in[0, t_0]$. We prove the conclusion in this step by contradiction. To this end, assume that there exist $t_2>t_1\geq t_0$ such that
\begin{equation}\label{eq36}
  \left\{
    \begin{aligned}
      Z(t_1) \in & ~ \partial \Omega^0, \\
      Z(t)\in &  ~ \Omega^1, ~t\in[t_1,t_2],\\
       Z(t_2)\in & ~ \partial \Omega^1.
    \end{aligned}
  \right.
\end{equation}

Consider the scaled ESO estimation error $\eta=[\eta_1, \ldots,$ $\eta_{n+1}]^{\rm{T}}$ with
\begin{equation}\label{eq31}
  \eta_{i}=\frac{x_i-\widehat{x}_i}{\varepsilon^{n+1-i}}, ~1\leq i\leq n+1.
\end{equation}
By (\ref{eq1}), (\ref{eso}), and (\ref{eq27}), the dynamics of $\eta$ can be formulated as
\begin{equation}\label{eq37}
  \left\{ \begin{aligned}
  \varepsilon\dot{\eta}_i=&-l_i\eta_1+\eta_{i+1}, ~1\leq i\leq n-1, \\
  \varepsilon\dot{\eta}_n=&-l_n\eta_{1}+\eta_{n+1}+\vartheta_1(\cdot),\\
  \varepsilon\dot{\eta}_{n+1}=&-l_{n+1}\eta_{1}-\vartheta_2(\cdot)\eta_1+\varepsilon\vartheta_3(\cdot),
  \end{aligned}
  \right.
\end{equation}
where
\begin{align*}
&  \vartheta_1(x,\widehat{x},\widehat{x}_{n+1},\widehat{\Theta}_c)=  f_0(x)-f_0(\widehat{x})\\
&  \qquad +\left(g_0(x)-g_0(\widehat{x})\right) \left[\widehat{u}_0(\overline{x},\widehat{\Theta}_c)-\frac{\overline{x}_{n+1}}{g_0(\overline{x})}\right], \\
&   \vartheta_2(x,z,\omega)= \frac{g(x,z,\omega)-g_0(x)}{g_0(\overline{x})}s'(\widehat{x}_{n+1}/M_{n+1}),  \\
&  \vartheta_3(x,z,\omega,\widehat{x},\widehat{x}_{n+1},\widehat{\Theta}_c)= \dot{f}(x,z,\omega)-\dot{f}_0(x) \\
& \qquad +(\dot{g}(x,z,\omega)-\dot{g}_0(x))\left[\widehat{u}_0(\overline{x},\widehat{\Theta}_c)-\frac{\overline{x}_{n+1}}{g_0(\overline{x})}\right] \\
   & \qquad +(g(x,z,\omega)-g_0(x))\dot{\widehat{u}}_0(\overline{x},\widehat{\Theta}_c) \\
   & \qquad +(g(x,z,\omega)-g_0(x))\frac{\overline{x}_{n+1}\dot{\overline{x}}}{g_0^2(\overline{x})}. \end{align*}
The equations of $\eta$ can be rewritten into the following compact form:
\begin{equation}\label{eq38}
\dot{\eta}=\frac{1}{\varepsilon}E\eta-F_1\vartheta_2(\cdot)\eta_1+\left[F_2\frac{\vartheta_1(\cdot)}{\varepsilon}+F_1\vartheta_3(\cdot)\right],
\end{equation}
where $F_1=[0 ~B^{\rm{T}}]^{\rm{T}}$, $F_2=[B^{\rm{T}} ~0]^{\rm{T}}$. The matrix $E$ is Hurwitz by  design. To show the convergence of $\eta$, we first consider the following system which is simplified from (\ref{eq38}):
\begin{equation}\label{eq39}
\dot{\eta}=\frac{1}{\varepsilon}E\eta-F_1\vartheta_2(\cdot)\eta_1.
\end{equation}
The system above can be regarded as a negative feedback connection of the time-varying gain $\theta_2(\cdot)$ and the transfer function
\begin{equation}\label{eq40}
  G(\varepsilon s)=\frac{l_{n+1}}{(\varepsilon s)^{n+1}+l_1(\varepsilon s)^{n}+\cdots+\l_{n+1}}.
\end{equation}
Note that $|\vartheta_2(\cdot)|\leq k_g|s'(\widehat{x}_{n+1}/M_{n+1})|\leq k_g$. By Assumption A3 and (\ref{eq40}), one has $k_g\|G\|_{\infty}<1$. It then follows from the circle criterion \cite{Khalil-2002} that the origin of (\ref{eq39}) is globally exponentially stable. Therefore, applying a loop transformation to (\ref{eq39}) and using the Kalman-Yakubovich-Popov lemma \cite{Khalil-2002}, there exists a quadratic Lyapunov function $W(\eta)=\eta^{\rm{T}}P\eta$ for system (\ref{eq39}) satisfying $\dot{W}(\eta)\leq -\frac{\alpha}{\varepsilon}W(\eta)$, where $P$ is a  positive definite matrix and $\alpha$ is an $\varepsilon$-independent positive constant. Utilizing $W(\eta)$ as a Lyapunov function candidate for system (\ref{eq38}), and calculating its derivative yields
\begin{equation}\label{eq41}
  \dot{W}(\eta)\leq -\frac{\alpha}{\varepsilon}W(\eta)+2\overline{\lambda}\|\eta\|\left|F_2\frac{\vartheta_1(\cdot)}{\varepsilon}+F_1\vartheta_3(\cdot)\right|,
\end{equation}
where $\overline{\lambda}=\lambda_{\max}(P)$. Note that 1) the functions $f_0(\cdot)$, $g_0(\cdot)$, and $s(\cdot)$ are continuously differentiable with locally Lipschitz derivatives and globally bounded; 2) by (\ref{eq36}), $x$, $\widehat{\Theta}_v$, and $\widehat{\Theta}_c$ are all bounded in the time interval $[0, t_2]$; 3) by (\ref{eq31}), $|\widehat{x}_i|\leq |x_i|+\varepsilon^{n+1-i}|\eta_i|$, $1\leq i\leq n+1$. Based on these observations, we can conclude that $\frac{\vartheta_1(\cdot)}{\varepsilon}$ and $\vartheta_3(\cdot)$ are locally Lipschitz, and bounded from above by affine in $\|\eta\|$ functions,  uniformly in $\varepsilon$. Here, we take the term $f_0(x)-f_0(\widehat{x})$ as an example to further clarify this point. According to the Hadamard's lemma \cite{Nest-2006}, one has
\begin{align*}
 & \quad \frac{1}{\varepsilon}\left[f_0(x)-f_0(\widehat{x})\right]=\frac{1}{\varepsilon}(x-\widehat{x})\int_0^1\frac{\partial f}{\partial x}(\widehat{x}+\lambda(x-\widehat{x}))\textrm{d}\lambda \\
 &  =\left[\varepsilon^{n-1}\eta_1,\ldots, \varepsilon\eta_{n-1},\eta_n\right]^{\rm{T}}\int_0^1\frac{\partial f}{\partial x}(\widehat{x}+\lambda(x-\widehat{x}))\textrm{d}\lambda.
\end{align*}
Therefore, from (\ref{eq41}), one has that there exists an $\varepsilon_1>0$ such that for any $\varepsilon\in(0,\varepsilon_1)$ and $\tilde{t}_1\in(t_0/2,t_0)$,
\begin{equation}\label{eq42}
  \|\eta\|=O(\varepsilon), ~\forall t\in [\tilde{t}_1,t_2].
\end{equation}
By the convergence of the ESO, and the definition of $x_{n+1}$, the bounds to saturate the output of the observer are selected to satisfy
\begin{align}
\label{eq52}  M_i>&\sup_{x\in\mathcal{X}}|x_i|, ~1\leq i\leq n, \\
  M_{n+1}>& \sup_{(x,z,\omega)\in \Xi, \widehat{\Theta}_c\in\Omega_c^1}\left|\frac{\Delta f(x,z,\omega)}{g(x,z,\omega)}\right. \nonumber \\
\label{eq53}  & \qquad \qquad\qquad ~ \left.+\frac{\Delta g(x,z,\omega)\widehat{u}_0(x,\widehat{\Theta}_c)}{g(x,z,\omega)}\right|.
\end{align}
It follows that the saturation elements $M_is(\widehat{x}_i/M_i)$, $1\leq i\leq n+1$,  all work in the linear zone in the time interval $[\tilde{t}_1,t_2]$, i.e., $\overline{x}_i=\widehat{x}_i$.

In the following, we consider the dynamics of $Z$. To facilitate the subsequent analysis, the NN-based approximations $\widehat{V}_x(\overline{x},\widehat{\Theta}_v)$ and $\widehat{u}_0(\overline{x},\widehat{\Theta}_c)$ are expressed in terms of $\widetilde{\Theta}_v$, $\widetilde{\Theta}_c$, $\kappa(\widehat{x})$, and $\kappa_x(\widehat{x})$ as
\begin{align}
  \widehat{V}_x\left(\overline{x},\widehat{\Theta}_v\right)= & \Theta^{\rm{T}} \phi_x(\widehat{x})+\left(\widehat{\Theta}_v-\Theta\right)^{\rm{T}}\phi_x(\widehat{x}) \nonumber \\
\label{eq34}  =& V^*_x(\widehat{x})+\widetilde{\Theta}^{\rm{T}}_v\phi_x(\widehat{x})-\kappa_x(\widehat{x}),\\
    \widehat{u}_0\left(\overline{x},\widehat{\Theta}_c\right)
  =& -\frac{1}{2}R^{-1}g_0^{\rm{T}}(\widehat{x})B^{\rm{T}}\phi^{\rm{T}}_x(\widehat{x})\left(\Theta+\widetilde{\Theta}_c\right) \nonumber\\
  =& u_0^*(\widehat{x})-\frac{1}{2}R^{-1}g_0^{\rm{T}}(\widehat{x})B^{T}\phi_x^{\rm{T}}(\widehat{x})\widetilde{\Theta}_c \nonumber\\
\label{eq35}   &  +\frac{1}{2}R^{-1}g^{\rm{T}}_0(\widehat{x})B^{\rm{T}}\kappa^{\rm{T}}_x(\widehat{x}).
\end{align}
By (\ref{eq5}) and (\ref{eq50}), and inserting (\ref{eq34}) and (\ref{eq35}) into (\ref{eq7}), the instantaneous BE $\delta_t$ can be written as
\begin{align}\label{eq33}
  \delta_t=  & \mu^{\rm{T}}\widetilde{\Theta}_v-\kappa_x(\widehat{x})\left[A\widehat{x}+B\left(f_0(\widehat{x})+g_0(\widehat{x})\widehat{u}_0(\widehat{x},\widehat{\Theta}_c)\right)\right] \nonumber \\
  & +2u_0^*(\widehat{x})R\left[-\frac{1}{2}R^{-1}g_0^{\rm{T}}(\widehat{x})B^{\rm{T}}\phi_x^{\rm{T}}(\widehat{x})\widetilde{\Theta}_c\right. \nonumber \\
  & \left .+\frac{1}{2}R^{-1}g^{\rm{T}}_0(\widehat{x})B^{\rm{T}}\kappa^{\rm{T}}_x(\widehat{x})\right]\nonumber \\
  & +\frac{R}{4}\left[-R^{-1}g_0^{\rm{T}}(\widehat{x})B^{\rm{T}}\phi_x^{\rm{T}}(\widehat{x})\widetilde{\Theta}_c+R^{-1}g^{\rm{T}}_0(\widehat{x})B^{\rm{T}}\kappa^{\rm{T}}_x(\widehat{x})\right]^2\nonumber \\
 = & \mu^{\rm{T}}\widetilde{\Theta}_v+\frac{1}{4}\widetilde{\Theta}^{\rm{T}}_cG_{t}\widetilde{\Theta}_c+\Delta,
\end{align}
where
\begin{align*}
   \Delta \triangleq & -\kappa_x(\widehat{x})[A\widehat{x}+Bf_0(\widehat{x})]\nonumber \\
   & +\frac{1}{4}g_0^{\rm{T}}(\widehat{x})B^{\rm{T}}\kappa^{\rm{T}}_x(\widehat{x})R^{-1}\kappa_x(\widehat{x})Bg_0(\widehat{x})\nonumber \\
   & +\frac{1}{2}\Theta^{\rm{T}}g_0^{\rm{T}}(\widehat{x})B^{\rm{T}}R^{-1}Bg_0(\widehat{x})\phi_x(\widehat{x})\kappa^{\rm{T}}_x(\widehat{x}).
\end{align*}
Similarly, the  extrapolated BE $\delta_i$ can be expressed as
\begin{equation}\label{eq19}
  \delta_i=\mu_i^{\rm{T}} \widetilde{\Theta}_v+\frac{1}{4}\widetilde{\Theta}^{\rm{T}}_cG_{i}\widetilde{\Theta}_c+\Delta_i,
\end{equation}
where
\begin{align*}
 \Delta_i \triangleq & -\kappa_x(x^i)[Ax^i+Bf_0(x^i)] \nonumber \\
 & +\frac{1}{4}g_0^{\rm{T}}(x^i)B^{\rm{T}}\kappa^{\rm{T}}_x(x^i)R^{-1}\kappa_x(x^i)Bg_0(x^i) \nonumber \\
 & +\frac{1}{2}\Theta^{\rm{T}}g_0^{\rm{T}}(x^i)B^{\rm{T}}R^{-1}Bg_0(x^i)\phi_x(x^i)\kappa^{\rm{T}}_x(x^i). \end{align*}

By (\ref{eq35}), the derivative of $V^*(x)$ under the control (\ref{eq27}) is given by
\begin{align}\label{eq43}
  \dot{V}^*(x)= & V_x^*(x)[Ax+B(f(x,z,\omega)+g(x,z,\omega)u)] \nonumber \\
  = & V_x^*(x)[Ax+B(f_0(x)+g_0(x)u_0^*(x))] \nonumber \\
    & +V_x^*(x)B\left[\left(x_{n+1}-\frac{g_0(x)}{g_0(\widehat{x})}\widehat{x}_{n+1}\right)\right. \nonumber \\
    &  +g_0(\widehat{x})\left(\widehat{u}_0(\widehat{x},\widehat{\Theta}_c)-u_0^*(\widehat{x})\right)\nonumber \\
    & +g_0(\widehat{x})\left(u_0^*(\widehat{x})-u_0^*(x)\right)\Big].
\end{align}
By (\ref{eq5}) and (\ref{eq33}), and the fact that the functions $V_x^*$, $\phi_x$, $g_0$, and $u_0^*$, are locally Lipschitz, one has that for all $x\in \mathcal{X}$, the derivative of $V^*(x)$ is upper bounded by
\begin{equation}\label{eq44}
\dot{V}^*(x)\leq -\iota_1\|x\|^2+\iota_2\|\eta\|+\iota_3\|\widetilde{\Theta}_c\|+\iota_4\overline{\kappa},
\end{equation}
where $\iota_1=\lambda_{\min}(Q)$, and $\iota_2$ to $\iota_4$ are $\varepsilon$-independent positive constants.

Let $\mathcal{V}: \mathbb{R}^n\times\mathbb{R}^{L}\times \mathbb{R}^{L}\rightarrow \mathbb{R}_{\geq 0}$ be a positive definite continuously differentiable Lyapunov function candidate defined by
\begin{equation}\label{eq32}
  \mathcal{V}(Z)=V^*(x)+V_v(\widetilde{\Theta}_v)+V_c(\widetilde{\Theta}_c).
\end{equation}
It follows from (\ref{eq15})-(\ref{eq17}) that the derivative of $\mathcal{V}$ satisfies
\begin{align}\label{eq32}
  \dot{\mathcal{V}}= & \dot{V}^*(x)+ \widetilde{\Theta}_v^{\rm{T}}\left(-\lambda_{v1}\frac{\mu}{\rho}\delta_t-\frac{\lambda_{v2}}{N} \sum_{i=1}^{N}\frac{\mu_i}{\rho_i}\delta_i\right)\nonumber \\
  & +\widetilde{\Theta}_c^{\rm{T}}\left(-\lambda_{c1}(\widehat{\Theta}_c-\widehat{\Theta}_v)-\lambda_{c2}\widehat{\Theta}_c\right)\nonumber\\
  & +\widetilde{\Theta}_c^{\rm{T}}\left(\frac{\lambda_{v1}G_{t}^{\rm{T}}\widehat{\Theta}_c\mu^{\textrm{T}}}{4\rho} + \sum_{i=1}^{N}\frac{\lambda_{v2}G^{\rm{T}}_{i}\widehat{\Theta}_c\mu^{\rm{T}}_i}{4N\rho_i}\right)\widehat{\Theta}_v\nonumber \\
  & -\frac{1}{2}\widetilde{\Theta}_v^{\rm{T}}\Gamma^{-1}\left(\beta\Gamma-\lambda_{v1}\frac{\Gamma\mu\mu^{\rm{T}}\Gamma}{\rho^2}\right)\Gamma^{-1}\widetilde{\Theta}_v. \end{align}
Using (\ref{eq33}) and (\ref{eq19}),  one has
\begin{align}\label{eq45}
  \dot{\mathcal{V}}= & \dot{V}^*(x) -\lambda_{v1}\widetilde{\Theta}_v^{\rm{T}}\frac{\mu\mu^{\rm{T}}}{\rho}\widetilde{\Theta}_v-\widetilde{\Theta}_v^{\rm{T}}\lambda_{v1}\frac{\mu}{\rho}\Delta \nonumber  \\
   & -\frac{\lambda_{v2}}{N}\widetilde{\Theta}_v^{\rm{T}}\sum_{i=1}^N\frac{\mu_i\mu_i^{\rm{T}}}{\rho_i}\widetilde{\Theta}_v -\widetilde{\Theta}_v^{\rm{T}}\frac{\lambda_{v2}}{N}\sum_{i=1}^N\frac{\mu_i}{\rho_i}\Delta_i \nonumber \\
   & -(\lambda_{c1}+\lambda_{c2})\|\widetilde{\Theta}_c\|^2+\lambda_{c1}\widetilde{\Theta}_c^{\rm{T}}\widetilde{\Theta}_v-\lambda_{c2}\widetilde{\Theta}_c^{\rm{T}}\Theta \nonumber \\
   & +\left(\|\Theta\|^2\widetilde{\Theta}_c^{\rm{T}}+\|\widetilde{\Theta}_c\|^2\Theta^{\rm{T}}+\widetilde{\Theta}_c^{\rm{T}}\Theta\widetilde{\Theta}_v\right)\nonumber \\
   & \times \left(\frac{\lambda_{v1}G_{t}^{\rm{T}}\mu^{\textrm{T}}}{4\rho} + \sum_{i=1}^{N}\frac{\lambda_{v2}G^{\rm{T}}_{i}\mu^{\rm{T}}_i}{4N\rho_i}\right)\nonumber \\
   & -\frac{1}{2}\beta\widetilde{\Theta}_v^{\rm{T}}\Gamma^{-1}\widetilde{\Theta}_v+\lambda_{v1}\widetilde{\Theta}_v^{\rm{T}}\frac{\mu\mu^{\rm{T}}}{\rho^2}\widetilde{\Theta}_v.
\end{align}
By some straightforward manipulations using (\ref{eq54}) and the Young's inequality, one can obtain $\|\frac{\mu}{\rho}\|\leq \frac{1}{2\sqrt{\gamma \varsigma_0}}$, $\|\frac{\mu_i}{\rho_i}\|\leq \frac{1}{2\sqrt{\gamma \varsigma_0}}$, $\|\frac{\mu\mu^{\rm{T}}}{\rho}\|\leq \frac{1}{\gamma \varsigma_0}$, $\|\frac{\mu_i\mu_i^{\rm{T}}}{\rho_i}\|\leq \frac{1}{\gamma \varsigma_0}$, $\|\frac{\mu_i\mu_i^{\rm{T}}}{\rho^2}\|\leq \frac{1}{4\gamma \varsigma_0}$. These together with (\ref{eq44}) and Assumption A4 yield
\begin{align}\label{eq46}
 \dot{\mathcal{V}} \leq & -\iota_1\|x\|^2+\iota_3\|\widetilde{\Theta}_c\|+\frac{\lambda_{v1}}{\gamma\varsigma_0}\|\widetilde{\Theta}_v\|^2  \nonumber \\
   & -\lambda_{v2}c\|\widetilde{\Theta}_v\|^2-(\lambda_{c1}+\lambda_{c2})\|\widetilde{\Theta}_c\|^2\nonumber \\
   &  +\frac{\lambda_{c1}}{2}(\|\widetilde{\Theta}_c\|^2+\|\widetilde{\Theta}_v\|^2)+\lambda_{c2}\|\Theta\| \|\widetilde{\Theta}_c\| \nonumber \\
   & +\chi_1\|\Theta\|^2\|\widetilde{\Theta}_c\|+\chi_1\|\Theta\|\|\widetilde{\Theta}_c\|^2 \nonumber \\
   & +\frac{\chi_1\|\Theta\|}{2}(\|\widetilde{\Theta}_c\|^2+\|\widetilde{\Theta}_v\|^2) \nonumber \\
   & -\frac{1}{2}\beta\varsigma_1^{-1}\|\widetilde{\Theta}_v\|^2+\frac{\lambda_{v1}}{4\gamma\varsigma_0}\|\widetilde{\Theta}_v\|^2 \nonumber \\
   & +\iota_2\|\eta\|+\iota_4\overline{\kappa}+\frac{\lambda_{v1}}{2\sqrt{\gamma\varsigma_0}}|\Delta|\|\widetilde{\Theta}_v\|\nonumber \\
   & +\frac{\lambda_{v2}}{2\sqrt{\gamma\varsigma_0}}\|\widetilde{\Theta}_v\| \sum_{i=1}^{N}|\Delta_i| \nonumber \\
 = & -\iota_1\|x\|^2+\iota_5\|\widetilde{\Theta}_v\|^2+\iota_6\|\widetilde{\Theta}\|_c^2 +\iota_7\|\widetilde{\Theta}_c\|+\iota_8,
\end{align}
where
\begin{align*}
 \iota_5 \triangleq & -\lambda_{v2}c-\frac{\beta}{2\varsigma_1}+\frac{\lambda_{c1}}{2}+\frac{\chi_1\|\Theta\|}{2}+\frac{5\lambda_{v1}}{4\gamma\varsigma_0}, \\
 \iota_6 \triangleq &-\frac{\lambda_{c1}}{2}-\lambda_{c2}+\frac{3}{2}\chi_1\|\Theta\|, \\
 \iota_7 \triangleq & \iota_3+\lambda_{c2}\|\Theta\|+\chi_1\|\Theta\|^2, \\
 \iota_8 \triangleq & \iota_2\|\eta\|+\iota_4\overline{\kappa}+\frac{\lambda_{v1}}{2\sqrt{\gamma\varsigma_0}}|\Delta|\|\widetilde{\Theta}_v\|\\
 & +\frac{\lambda_{v2}}{2\sqrt{\gamma\varsigma_0}}\|\widetilde{\Theta}_v\|\sum_{i=1}^{N}|\Delta_i|,\\
 \chi_1 \triangleq & \frac{\lambda_{v1}}{8\sqrt{\gamma \varsigma_0}}\sup_{\overline{x}\in\mathcal{X}}\|G_{t}\|+\frac{\lambda_{v2}}{8\sqrt{\gamma \varsigma_0}}\max_{1\leq i \leq N}\|G_{i}\|.
\end{align*}
Sufficient conditions for the gains $\lambda_{c1}$, $\lambda_{v1}$, $\lambda_{v2}$, and $\beta$  are then given by
\begin{equation}\label{eq47}
 \beta >\varsigma_1 \chi_1\|\Theta\|+\varsigma_1 \lambda_{c1}, ~\lambda_{v2} > \frac{5\lambda_{v1}}{4c\gamma\varsigma_0}, ~\lambda_{c1} >3\chi_1\|\Theta\|.
\end{equation}
Provided the gains are selected satisfying (\ref{eq47}), the ESO is convergent in the sense of (\ref{eq42}),
the upper bound of the approximation error $\kappa_x$ can be made sufficiently small by increasing the number of NN neurons, and (\ref{eq36}), the derivative of $\mathcal{V}$ satisfies
\begin{equation}\label{eq48}
  \dot{\mathcal{V}}(Z(t))\leq 0,  ~t\in[t_1,t_2].
\end{equation}
This contradicts (\ref{eq36}). Thus the statement in Step 1) holds.

\emph{Step 2): There exists $\varepsilon^*\in(0, \varepsilon^{\dag}]$ such that for any $\varepsilon\in(0,\varepsilon^*)$, $\|\eta(t)\|=O(\varepsilon)$, $\forall t\in[T,\infty)$, $T>0$, and $Z(t)$ is uniformly ultimately bounded.}

Since $Z(t)\in\Omega^1$, $\forall t\in[0, \infty)$, by (\ref{eq41}), one has that for any $\varepsilon\in(0,\varepsilon^*)$ and $T>0$, $\|\eta(t)\|=O(\varepsilon)$, $\forall t\in[T,\infty)$. What is more, (\ref{eq46}) holds for $t\in[T,\infty)$. Consequently, selecting the parameters according to (\ref{eq47}), and by(\ref{eq46}), one can conclude the uniformly ultimately boundedness of $Z$. Also note that the system state $x$ converges to the neighbourhood of origin, and the estimate weight $\widehat{\Theta}_c$ approximates the ideal weight $\Theta$ if $\iota_7/|\iota_6|$ (i.e., $\lambda_{c2}/\lambda_{c1}$) is sufficiently small. This completes the proof of Theorem 1. \IEEEQED

\emph{Remark 5:} In this paper, similar to \cite{Luo-2014,Kama-2016a,Jiang-2014}, we use two weights $\widehat{\Theta}_v$ and $\widehat{\Theta}_c$ to estimate the same ideal weight $\Theta$. The use of two weights rather than one weight is mainly motivated by the stability analysis. Note that in (\ref{eq33}) and (\ref{eq19}), benefited from the use of two weights, the BEs are linear with respect to the weight estimation error $\widetilde{\Theta}_v$, which enables us to develop a least-square update law for $\widehat{\Theta}_v$ based on the BEs. The update law for $\widehat{\Theta}_c$ guarantees that $\widehat{\Theta}_c\rightarrow \widehat{\Theta}_v$. $\IEEEQED$

\begin{figure}
  \centering
  \includegraphics[width=0.45\textwidth,bb=10 20 480 250, clip]{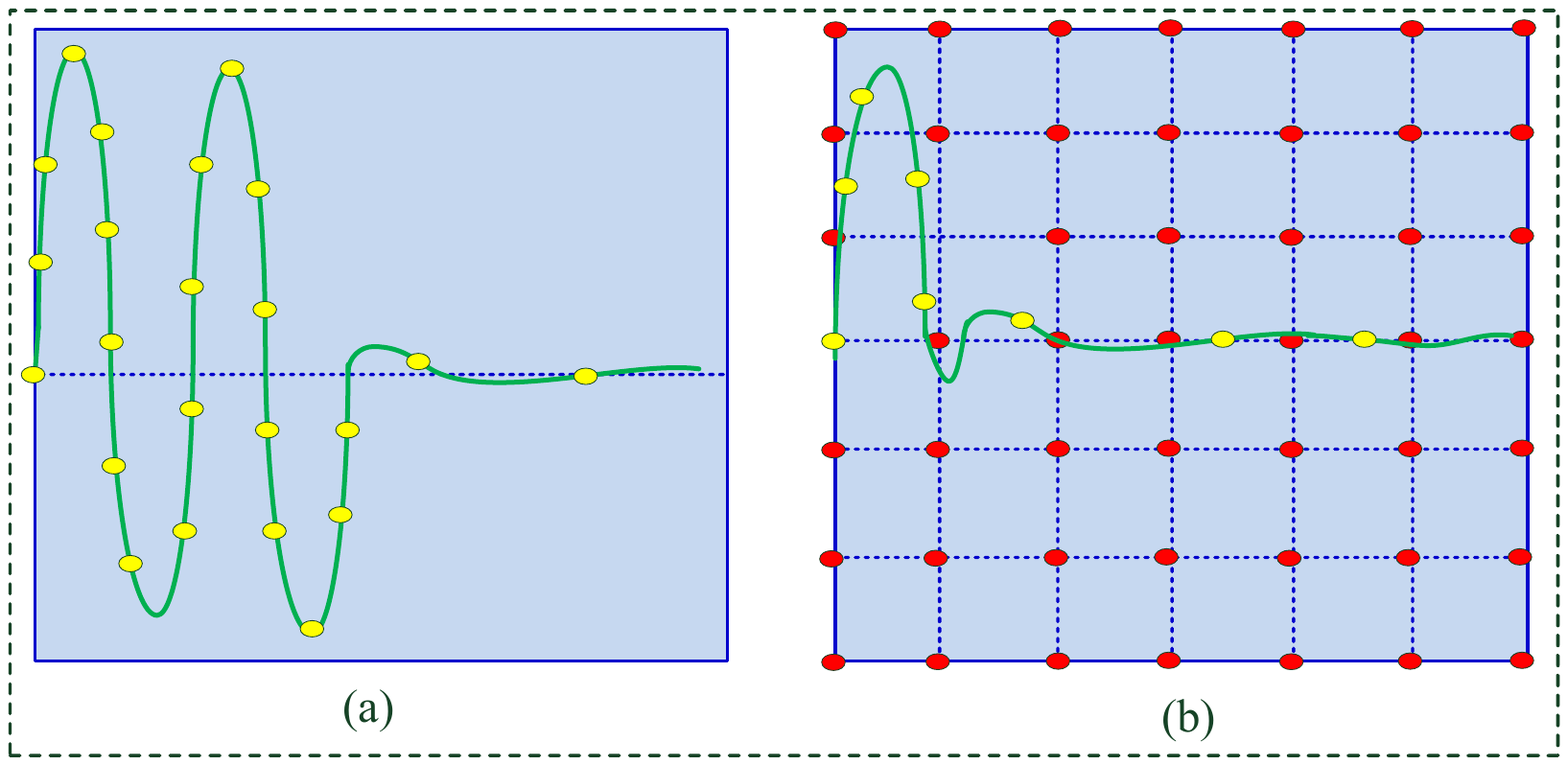}
   \caption{Illustration of the differences between the mechanisms of the standard CL update law and the update law in this paper: (a) standard CL update law; (b) update law in this paper.  The green lines represent the system trajectories. The yellow dots represent the data points along the trajectory. The red dots represent the selected data points to extrapolate the BE. }\label{data_points}
\end{figure}

\emph{Remark 6:} We should point out that the update law in this paper is different from the standard CL update law in existing wroks such as \cite{Zhang-2018,Xue-2020,Yang-2019b,Yang-2020a,Vam-2015,Chowdhary-2010}. As depicted in Fig. \ref{data_points}a, for a standard CL update law, the recorded data can be only selected along the system state trajectory. In this case, to guarantee that the system state visits sufficient number of points in the domain of operation, a probing signal is generally required to excite the system \cite{Vam-2015}. In this paper, with the known nominal model of the system, the BE can be extrapolated to any desired data points, and the system trajectory doesn’t need to really reach to these data points (see Fig. \ref{data_points}b). Benefited from this, our developed approach doesn't need a probing signal, and hence it is able to achieve better transient performance. \IEEEQED

\emph{Remark 7:} According to Theorem 1, $\widehat{x}\rightarrow x$ and the actor NN weight $\widehat{\Theta}_c$ approximates the ideal weight $\Theta$. Therefore, the control policy given by $\widehat{u}_0(\overline{x},\widehat{\Theta}_c)$ approximates the ideal optimal control policy $u_0^*(x)$. In this sense and following also \cite{Lewis-2009,Lewis-2018,Liu-2017}, $\widehat{u}_0(\overline{x},\widehat{\Theta}_c)$ is referred to as an \emph{approximate optimal control}. Also note that, similar to \cite{Kama-2016a,Kama-2016b,Kama-2013}, the convergence results in this paper are practical. Thus the cost defined by (\ref{eq3}), when evaluated along the approximate optimal trajectory, will be infinite. This, however, will not cause any theoretical concern.
Note that the partial derivative $\widehat{V}_x$ leveraged in the BEs (\ref{eq7})-(\ref{eq49}) to update the NN weights is bounded for all $t\in[0,\infty)$. Actually, the design and analysis of the RL-based controller doesn't need the infinite horizon cost along the approximate optimal trajectory to be finite, it only needs the cost, when evaluated along the ideal optimal trajectory to be finite. To avoid infinite cost, one could use a discounted cost formulation with some additional conditions to guarantee stability \cite{Gait-2015}. \IEEEQED

\emph{Remark 8:} There are several groups of parameters that need to be selected for the developed ESO-based reinforcement learning and disturbance rejection scheme.
\begin{itemize}
  \item  For the ESO, the parameters to be selected include the observer gains $L$ and $\varepsilon$, the saturation bounds for the observer output ($M_i$, $1\leq i\leq n+1$), and the saturation bounds for the functions $f_0(\cdot)$ and $g_0(\cdot)$ (denoted by $M_f$ and $M_g$, respectively). The observer gain $L$ can be selected by a pole placement method such that the matrix $E$ is Hurwitz. Theoretically, $\varepsilon$ can be selected arbitrarily small to achieve more accurate estimation. However, in practice, the lower-bound of $\varepsilon$ is limited due to noise and sampling constraints \cite{Khalil-2008,Ran-2017a}. The bounds $M_i$, $1\leq i\leq n+1$, are selected according to (\ref{eq52})-(\ref{eq53}). Note that the calculation of $M_{n+1}$ is generally not straightforward, and one might end up with a conservative bound. From a practice viewpoint, the values of $\varepsilon$ and $M_{n+1}$ can be decided by a simple trial and error procedure, based on the obtained performance. Our numerical experience and many previous ESO results (see, e.g., \citep{Khalil-2008,Ran-2017a}) indicate that it is generally very easy to select a group of satisfactory $\varepsilon$ and $M_{n+1}$. Finally, the bounds $M_f$ and $M_g$ are selected such that $M_f\geq \sup_{x\in\mathcal{X}}|f_0(x)|$ and $M_g\geq \sup_{x\in\mathcal{X}}|g_0(x)|$, respectively.
  \item For the simulation of experience based RL, the parameters to be selected include the positive adaptation gains $\lambda_{c1}$, $\lambda_{c2}$, $\lambda_{v1}$, $\lambda_{v2}$, and $\beta$. The sufficient conditions for the adaptation gains based on the stability analysis are given by (\ref{eq47}).  Note that the inequalities in (\ref{eq47}) depend on unknown parameters $\Theta$, $\varsigma_0$, $\varsigma_1$, and $c$. Therefore, the selection of these adaptation gains also needs a simple trial and error procedure. What is more, $\lambda_{c2}/\lambda_{c1}$ is required to be sufficiently small.
   \item For the data set $\{x^i\in\mathbb{R}^n| i=1,\ldots, N\}$ to extrapolate the BE, it is required to satisfy Assumption A4. We mention that, similar to the probing signals used in previous RL algorithms \cite{Jiang-2014b,Gao-2019,Fan-2016,Wang-2018,Wang-2020},  the data set is also selected offline and then the RL algorithm is implemented online. It is difficult to provide a theoretical guarantee to make sure that the offline selected data points (and probing signals) achieve satisfactory online performance. 
       In practice, to fulfill Assumption A4, one can select $x^i$, $1\leq i\leq N$, on an $\underbrace{a\times a \cdots \times a}_{n}$ data grid which covers the interested compact set $\mathcal{X}\subseteq\mathbb{R}^n$, where $a$ is an appropriately large positive integer.
   \item For the basis function $\phi(x)$, and the initial NN weights $\widehat{\Theta}_v(0)$ and $\widehat{\Theta}_c(0)$, their selections are important for RL-based controllers \cite{Lewis-2009,Lewis-2018,Liu-2017}.  The basic principle of selecting the basis function is to let the real basis be contained in the selected basis to make sure that the NN representation error $\kappa(x)$ is small. In this paper, since the system nominal model (i.e., $f_0(x)$ and $g_0(x)$) is known to the designer, the selection of the basis function will be relatively easier than that for an unknown system. However, the  selection of a good basis for very general nonlinear systems is challenging and is still largely open in machine learning \cite{Kama-2013}. On the other hand, similar to \cite{Kama-2016a,Kama-2016b}, benefited from the mechanism that the BE can be extrapolated to any desired data points with a known system nominal model, the initial NN weights are not required to be admissible. However, in practice, better initial NN weights will lead to better closed-loop transient performance. 
       \IEEEQED
\end{itemize}

\emph{Remark 9:} Compared with the existing RL results for uncertain nonlinear systems \cite{Jiang-2014b,Gao-2019,Fan-2016,Xue-2020,Jiang-2017b,Wang-2018,Wang-2020,Zhang-2018,Yang-2019b}, the advantages of the proposed approach are threefold: 1) The uncertainties considered in this paper are more general and complex, since it involves the system states $x$ and $z$, the control input $u$, and the external disturbance $\omega$, and hence the analysis is much more challenging.  2) This paper inherits the idea of ADRC and handles the uncertainty in an ``observation+compensation'' scheme, which is fundamentally different from the ideas in \cite{Jiang-2014b,Gao-2019,Fan-2016,Xue-2020,Jiang-2017b,Wang-2018,Wang-2020,Zhang-2018,Yang-2019b}. 3) This paper provides a more practical solution to the RL design problem for uncertain nonlinear systems since our developed approach is output feedback based and does not require the probing signal. The approaches in \cite{Jiang-2014b,Gao-2019,Jiang-2017b,Wang-2018,Wang-2020,Zhang-2018,Yang-2019b,Fan-2016,Xue-2020} are full state-feedback based and in \cite{Jiang-2014b,Gao-2019,Fan-2016,Wang-2018,Wang-2020} require the probing signal. \IEEEQED


\begin{figure}[!t]
  \centering
  \includegraphics[width=0.30\textwidth,bb=15 10 435 365, clip]{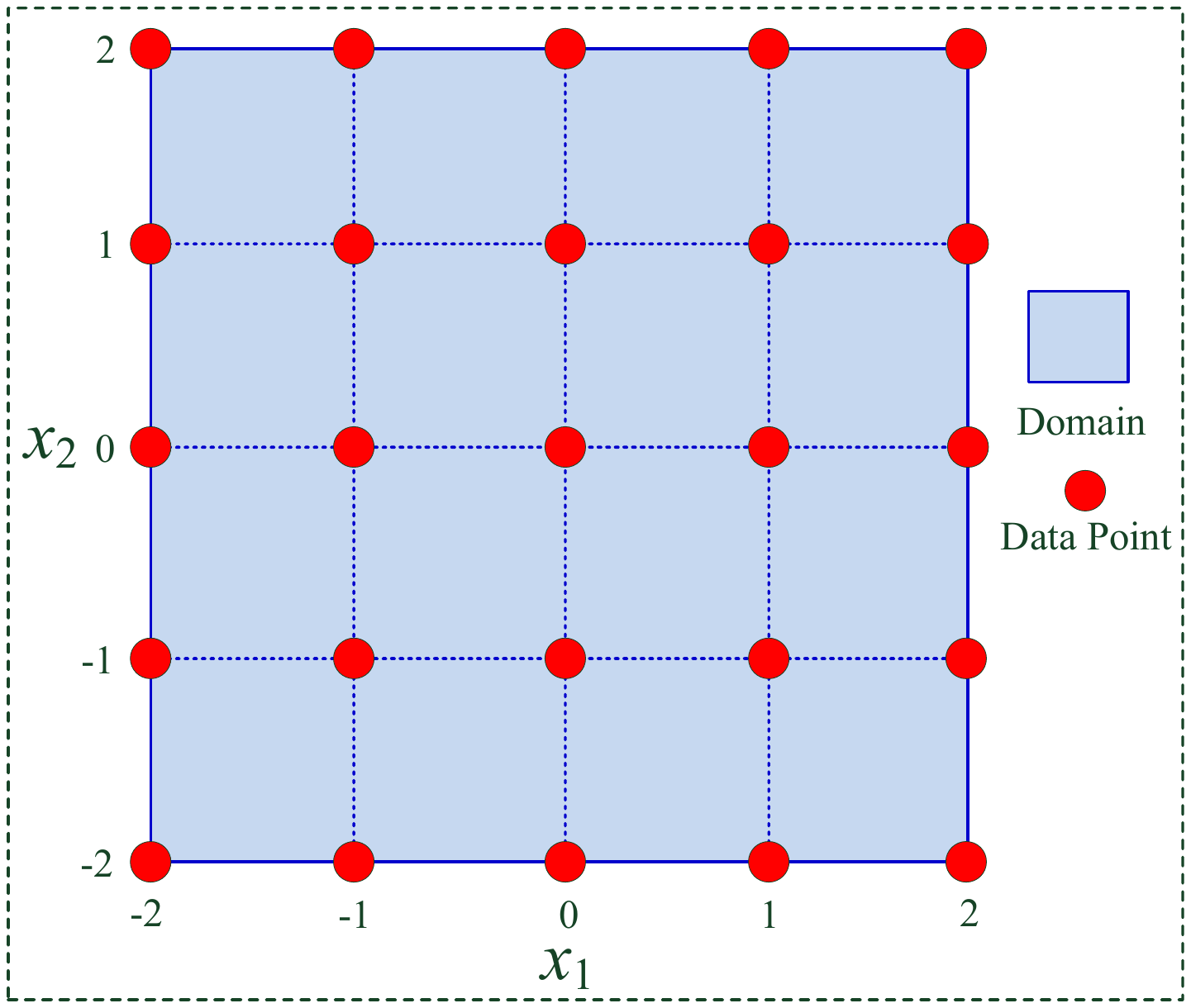}
   \caption{Illustration of the data points selected on a $5\times 5$ data grid.}\label{fig11}
\end{figure}

\begin{figure}
  \centering
  \includegraphics[width=0.40\textwidth,bb=15 10 315 382, clip]{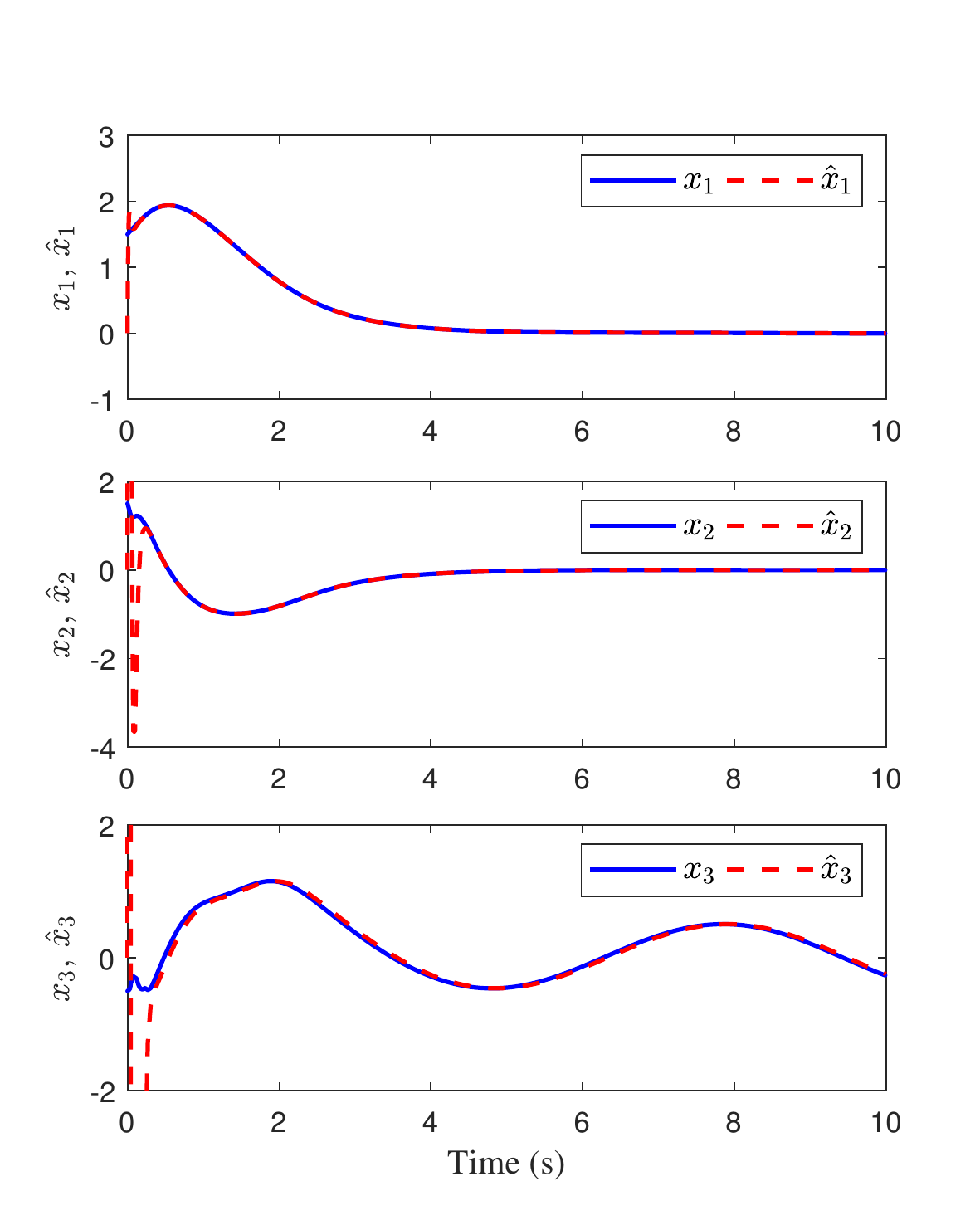}
   \caption{Trajectories of the system state and ESO output.}\label{fig3}
\end{figure}

\section{Examples}

This section presents two simulation examples to illustrate the effectiveness of our developed approach.

\begin{figure}
  \centering
  \includegraphics[width=0.40\textwidth,bb=15 0 315 188, clip]{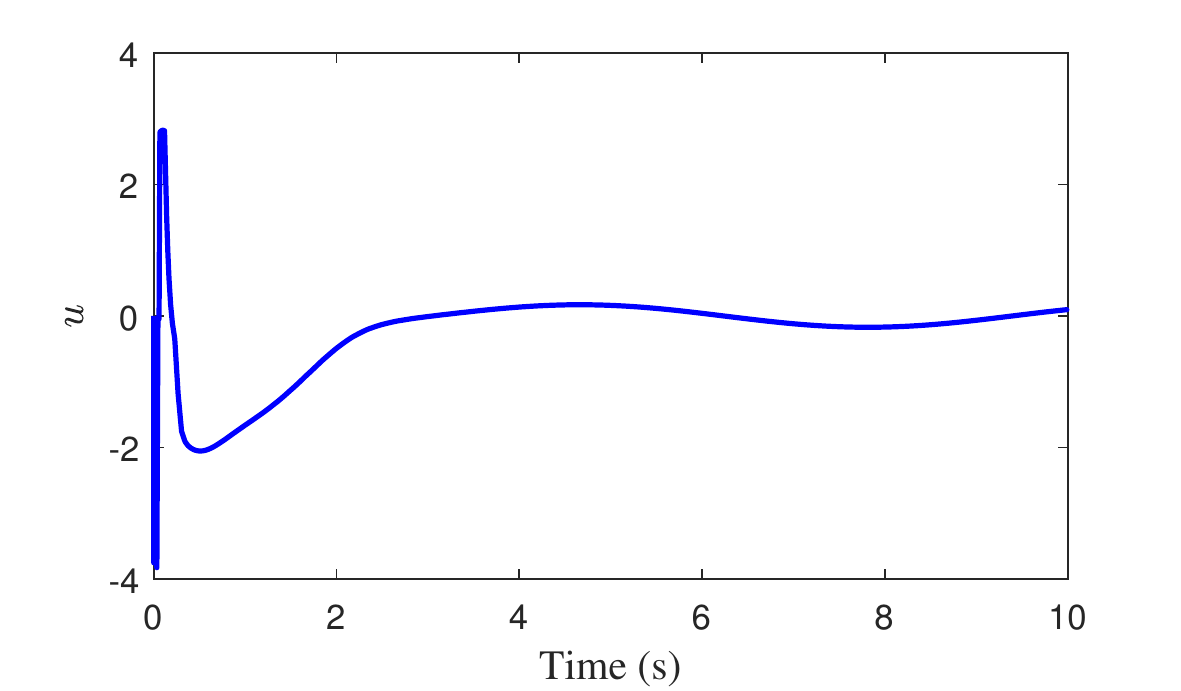}
  \caption{Trajectory of the control input $u$.}\label{fig4}
\end{figure}

\begin{figure}
  \centering
  \includegraphics[width=0.40\textwidth,bb=5 0 315 240, clip]{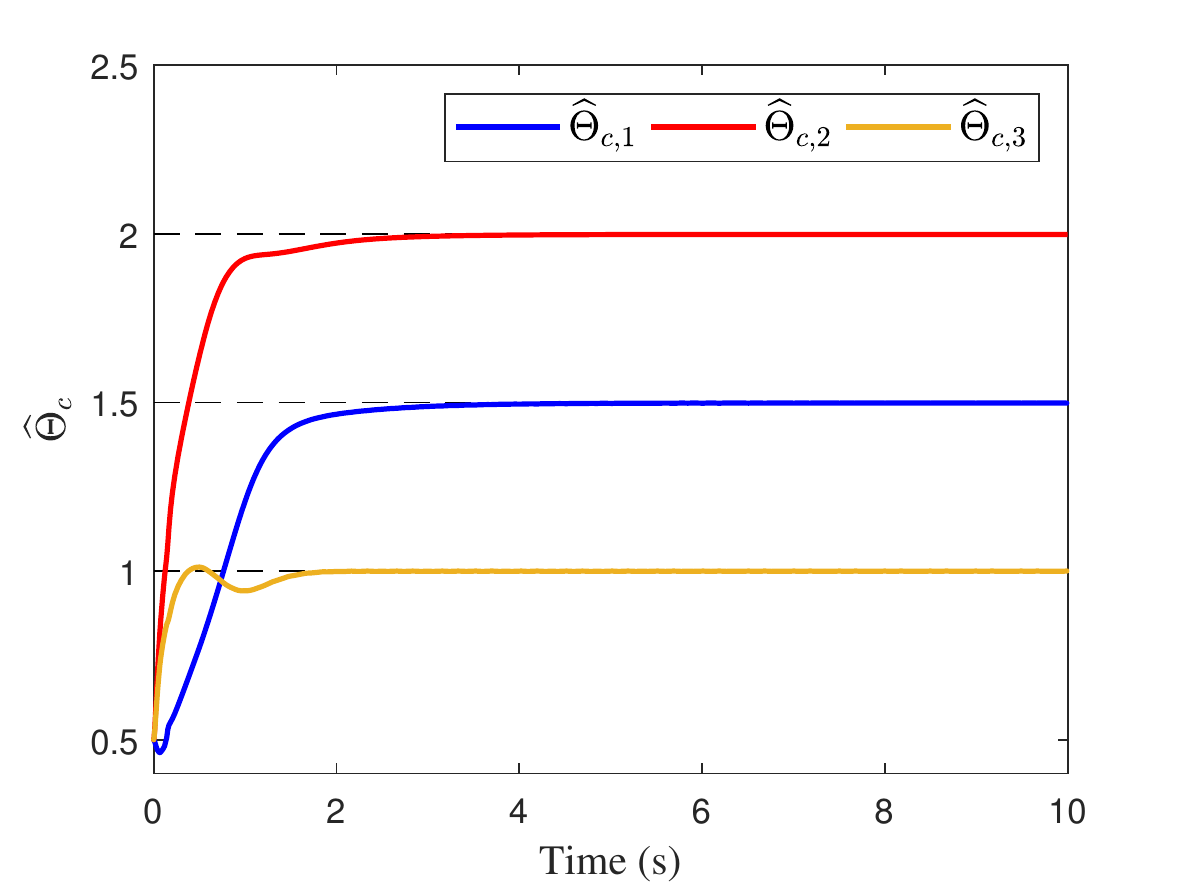}
  \caption{Trajectories of the actor NN weight.}\label{fig6}
\end{figure}

\subsection{Example 1}
Consider the following uncertain nonlinear system
\begin{equation}\label{eq55}
   \left\{
  \begin{aligned}
          \dot{z}=& \underbrace{-(x_2^2+\omega^2)z}_{f_z(x,z,\omega)}, \\
          \dot{x}_1=& x_2, \\
          \dot{x}_2=& \underbrace{-x_1-2.5x_2+\omega+z^2+0.5(x_1+x_2)(\cos(2x_1)+2)^2}_{f(x,z,\omega)} \\
          & +\underbrace{(\cos(2x_1)+2+\sin(x_1)\omega)}_{g(x,z,\omega)}u, \\
                y=& x_1,
        \end{aligned} \right.
\end{equation}
where the external disturbance is numerically set as $\omega=0.5\sin(t)$. The known nominal models of $f(x,z,\omega)$ and $g(x,z,\omega)$ are taken as
\begin{align*}
  f_0(x)= & -x_1-1.5x_2+0.5(x_1+x_2)(\cos(2x_1)+2)^2,  \\
  g_0(x)= & \cos(2x_1)+2.
\end{align*}
The cost functional for the  nominal system is given by (\ref{eq3}) with $Q=x^{\rm{T}}\overline{Q}x$ and $R=1$, where $\overline{Q}=\left[
                \begin{array}{cc}
                  2 & 1 \\
                  1 & 1 \\
                \end{array}
              \right]$
is positive definite. The  nominal models are selected above since the corresponding optimal control problem has an analytical solution, which is helpful for the simulation to verify the correctness of the developed approach. Specifically, according to the procedure in \cite{Nev-1996},
the optimal value function is $V^*(x)=1.5x_1^2+2x_1x_2+x_2^2$, and the optimal control policy is $u_0^*(x)=-(\cos(2x_1)+2)(x_1+x_2)$. The total uncertainty for system (\ref{eq55}) is denoted by
\begin{equation}\label{eq56}
  x_3=-x_2+\omega+z^2+\sin(x_1)\omega u.
\end{equation}
It can be verified that Assumptions A1 to A3 are satisfied.


\begin{figure}
  \centering
  \includegraphics[width=0.40\textwidth,bb=5 40 315 580, clip]{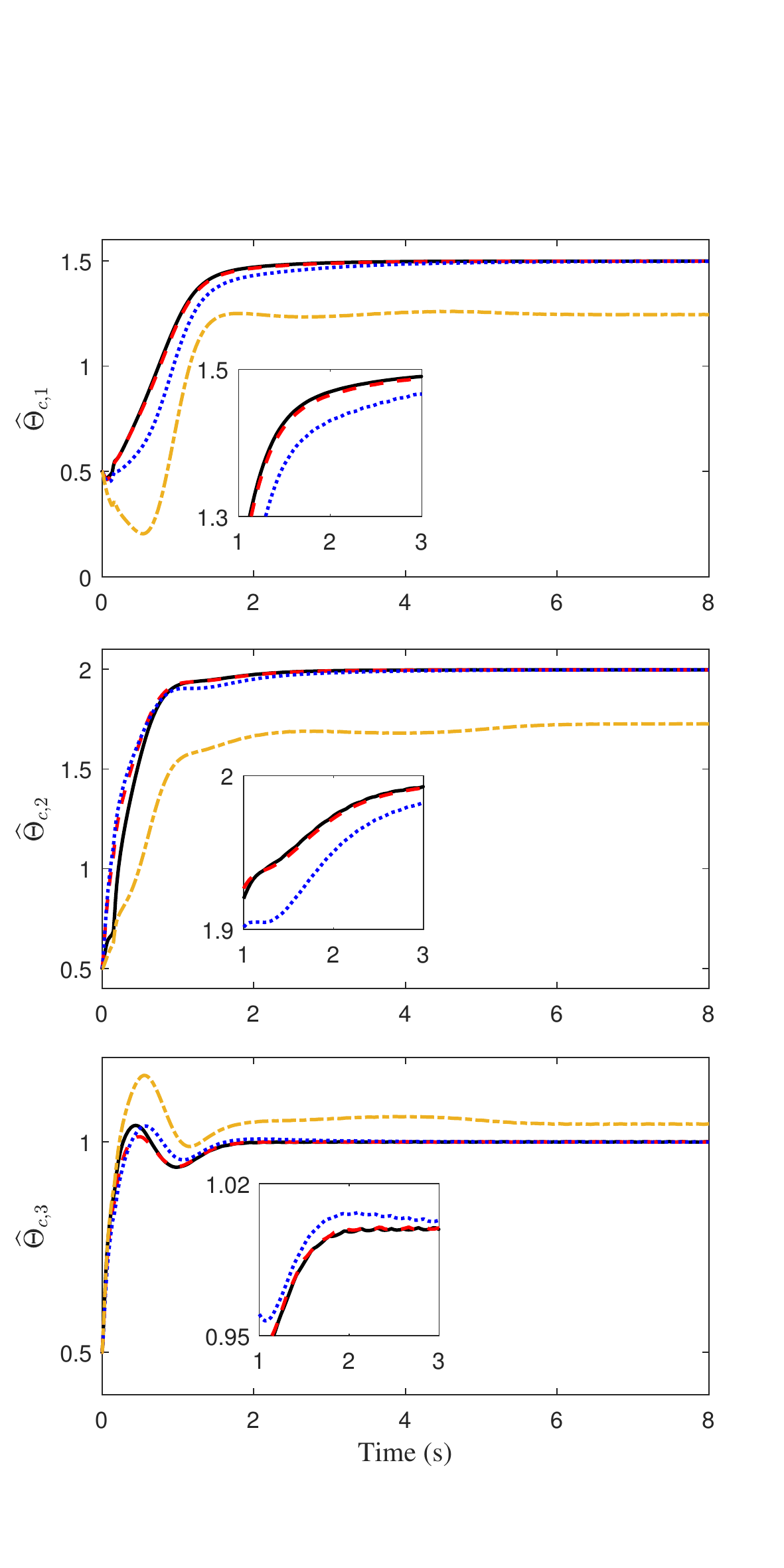}
  \caption{Trajectories of the actor NN weight with different data grids: $9\times 9$ data grid (black full line); $5\times 5$ data grid (red dash line); $3\times 3$ data grid (blue dot line); $2\times 2$ data grid (yellow dash-dot line).}\label{fig7}
\end{figure}

Consider the scenario that the compact set of interest is $\mathcal{X}=[-2 ~2]\times [-2 ~2]$. The ESO is designed with $L=[3 ~3 ~1]^{\rm{T}}$ and $\varepsilon=0.02$. The saturation bounds for the nominal functions $f_0(x)$ and $g_0(x)$ are selected as $M_f=6.5>\sup_{x\in\mathcal{X}}|f_0(x)|$ and $M_g=3>\sup_{x\in\mathcal{X}}|g_0(x)|$, respectively. The  saturation bounds for the output of the observer $\widehat{x}_1$, $\widehat{x}_2$, and $\widehat{x}_3$, are selected as $M_1=M_2=M_3=2$.

\emph{Simulation with known basis function:}
First, we simulate the simple case, i.e., the basis function is known. In this case, the basis function $\phi:\mathbb{R}^2\rightarrow \mathbb{R}^3$ is selected as $\phi(x)=[x_1^2 ~x_1x_2 ~x_2^2]^{\rm{T}}$. According to the analytical solution of $V^*(x)$, the ideal weight $\Theta=[1.5 ~2 ~1]^{\rm{T}}$. The data points to extrapolate the BE are selected on a $5\times 5$ data grid cover the domain $\mathcal{X}$ (see Fig. \ref{fig11}). The gains for the RL are selected as $\lambda_{v1}=1$, $\lambda_{v2}=5$,  $\lambda_{c1}=100$, $\lambda_{c2}=0.1$, $\gamma=0.5$, and $\beta=100$. The upper bound of $\Gamma$ is set as $\varsigma_1=2000$. Simulation is done with initial conditions $z(0)=1$, $x(0)=[1.5 ~1.5]^{\rm{T}}$, $[\widehat{x}_1(0) ~\widehat{x}_2(0) ~\widehat{x}_3(0)]^{\rm{T}}=[0 ~0 ~0]^{\rm{T}}$, $\widehat{\Theta}_v(0)=\widehat{\Theta}_c(0)=[0.5 ~0.5 ~0.5]^{\rm{T}}$, and $\Gamma(0)=\textrm{diag}\{100,100,100\}$. Note that the initial weight $[0.5 ~0.5 ~0.5]^{\rm{T}}$ is not an admissible control policy. Fig. \ref{fig3} shows the trajectories of the system state and ESO output. It can be observed that the system state $x$ is regulated to the origin, and the state $x$ and total uncertainty $x_3$ are well-estimated by the ESO. Fig. \ref{fig4} illustrates the trajectory of the control signal $u$ given by (\ref{eq27}). Fig. \ref{fig6} shows that the actor NN weight $\widehat{\Theta}_{c}$  converges to its real value.

\begin{figure}
  \centering
  \includegraphics[width=0.40\textwidth,bb=25 15 315 325, clip]{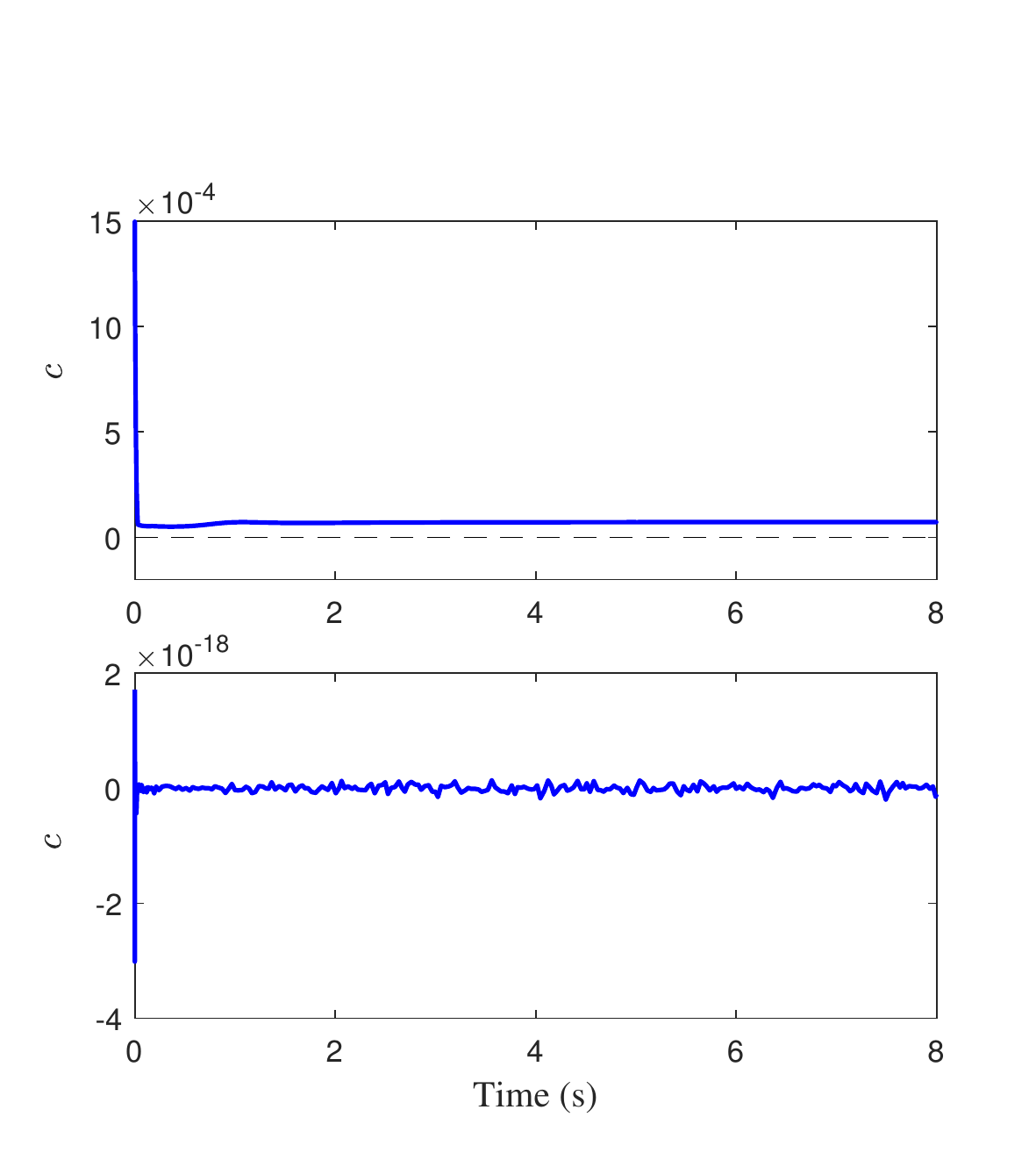}
  \caption{Trajectories of $\frac{1}{N}\left(\lambda_{\min}\left\{\sum_{i=1}^N\frac{\mu_i\mu_i^{\rm{T}}}{\rho_i}\right\}\right)$ with $3\times 3$ data grid (upper figure) and $2\times 2$ data grid (lower figure).}\label{fig8}
\end{figure}

\emph{Simulation with different data grids:} Next, we investigate the performance of the developed controller with different data grids. Note that one needs to select sufficiently more data points such that Assumption A4 is satisfied. However, too more data will increase the online computation burden. We consider the $9\times 9$, $5\times 5$, $3\times 3$, and $2\times 2$ data grids over the compact set $\mathcal{X}$. Fig. \ref{fig7} shows the trajectories of the actor NN weight with different data grids. From this figure, it can be observed that 1) more data points lead to better convergence performance; 2) the performance with $5\times 5$ data grid is comparable to $9\times 9$ data grid; 3) the actor NN weight cannot converge to the ideal value with $2\times 2$ data grid. Therefore, for this numerical example, $3\times 3$ data grid suffices to guarantee the convergence of the NN weights, and $5\times 5$ data grid is enough to obtain satisfactory  performance. Fig. \ref{fig8} depicts the trajectories of $\frac{1}{N}\left(\lambda_{\min}\left\{\sum_{i=1}^N\frac{\mu_i\mu_i^{\rm{T}}}{\rho_i}\right\}\right)$ with $3\times 3$ and $2\times 2$ data grids. One can see that Assumption A4 is satisfied with $3\times 3$ data grid, but failed with $2\times 2$ data grid.

\begin{figure}
  \centering
  \includegraphics[width=0.40\textwidth,bb=5 0 315 240, clip]{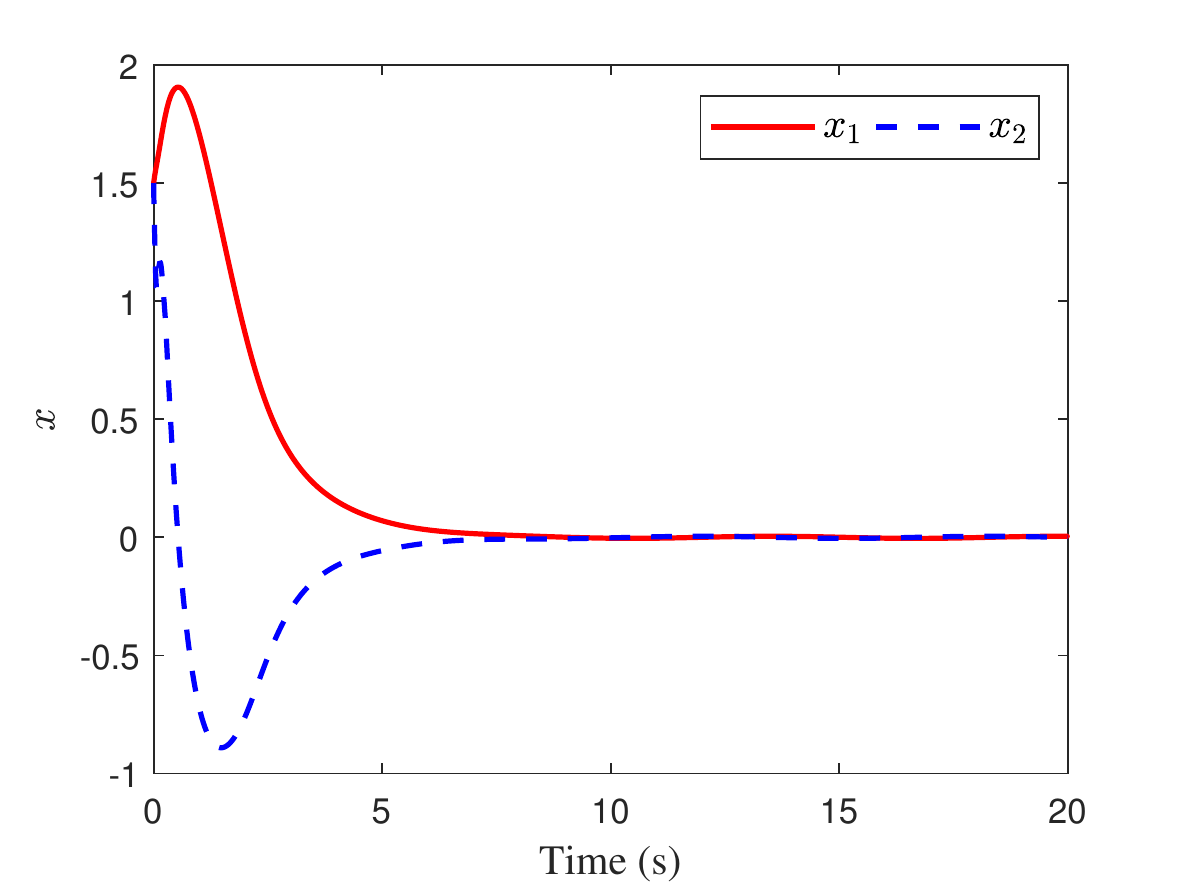}
  \caption{Trajectories of the system state with unknown basis function.}\label{fig9}
\end{figure}

\begin{figure}[!t]
  \centering
  \includegraphics[width=0.40\textwidth,bb=5 0 315 240, clip]{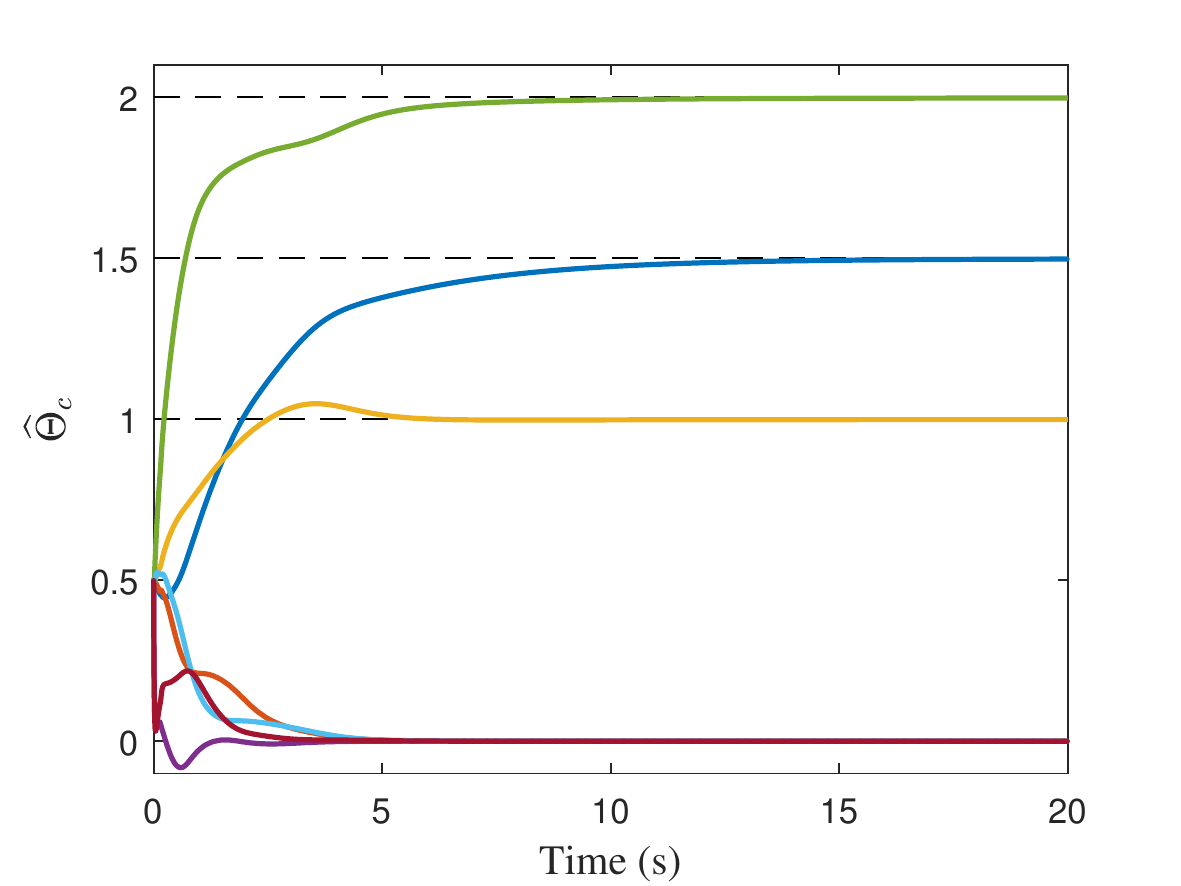}
  \caption{Trajectories of the actor NN weight with unknown basis function.}\label{fig10}
\end{figure}

\emph{Simulation with unknown basis function:}  Finally, we simulate the case that the basis function is unknown.  For this example, we select the basis function as $\phi(x)=[x_1^2 ~x_1^3 ~x_2^2 ~x_2^3 ~x_1x_2 ~x_1x_2^2 ~x_1^2x_2]^{\rm{T}}$. With this selection, the ideal value of the critical and actor NN weights is $\Theta=[1.5 ~0 ~1 ~0 ~2 ~0 ~0]^{\rm{T}}$. The initial conditions of $\widehat{\Theta}_v$, $\widehat{\Theta}_c$, and $\Gamma$, are set as $\widehat{\Theta}_v(0)=\widehat{\Theta}_c(0)=[0.5 ~0.5 ~0.5 ~0.5 ~0.5 ~0.5 ~0.5]^{\rm{T}}$, and $\Gamma(0)=\textrm{diag}\{100, 100, 100, 100, 100, 100, 100\}$. All other settings are the same with the simulation with known basis function. Trajectories of the system state $x$ and the actor NN weight $\widehat{\Theta}_c$ are illustrated in Figs. \ref{fig9} and \ref{fig10}, respectively. It can be observed that the system state converges to the origin while the actor NN weight converges to the ideal value. However, as expected, the transient period is longer than the known basis function case.

\begin{figure}
  \centering
  \includegraphics[width=0.40\textwidth,bb=10 25 315 475, clip]{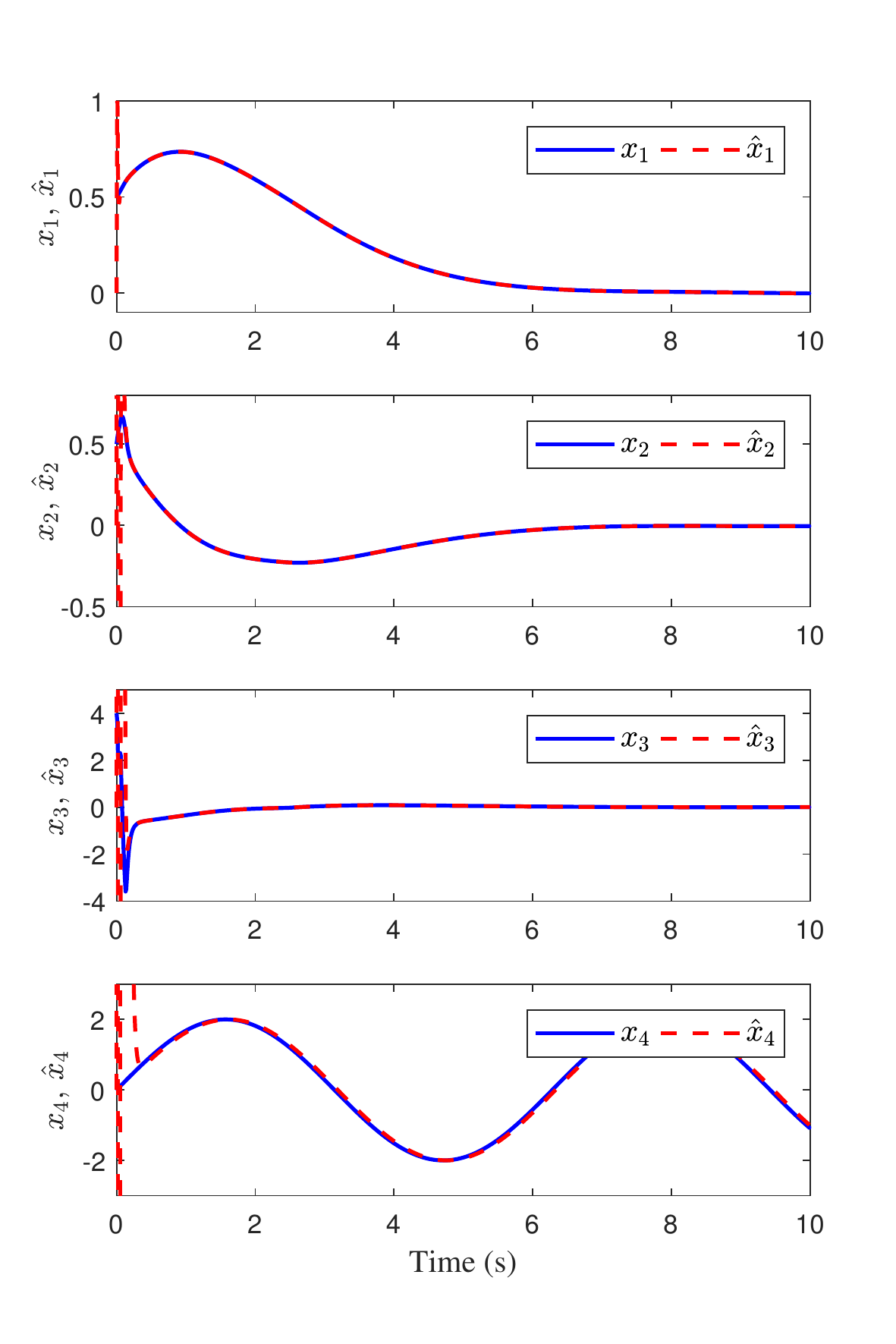}
   \caption{Trajectories of the system state and ESO output (Example 2).}\label{fig12}
\end{figure}

\subsection{Example 2}

In this example, we consider a class of third-order kinematic models described by  \cite{Bian-2019}
\begin{equation}\label{eq57}
 \left\{
  \begin{aligned}
          \dot{p}=& v, \\
          m\dot{v}=& f_a-bv, \\
                \tau \dot{f}_a=& u-f_a+\omega,
        \end{aligned} \right.
\end{equation}
where $p$, $v$, and $f_a$ denote the position, velocity, and actuator force, respectively; $u$ is the control input; $m$, $b$, and $\tau$ represent the mass, velocity constant, and time constant, respectively; and $\omega$ is the external disturbance. The above system represent a large class of physical systems, including human motor system, autonomous vehicle, power system, etc \cite{Bian-2019}. By letting $x_1=p$, $x_2=v$, and $x_3=\frac{f_a-bv}{m}$, the above system can be written as
\begin{equation}\label{eq58}
 \left\{
  \begin{aligned}
          \dot{x}_1=& x_2, \\
          \dot{x}_2=& x_3,\\
          \dot{x}_3=& -\frac{b}{\tau m}x_2-\left(\frac{1}{\tau}+\frac{b}{m}\right)x_3+\frac{1}{\tau m}u+ \frac{1}{\tau m}\omega.
        \end{aligned} \right.
\end{equation}

In the simulation, the values of the model parameters are given by  $m=1\textrm{kg}$, $b=1 \textrm{N} \cdot \textrm{s/m}$, and $\tau=0.1\textrm{s}$. Let $f_0(x)=-\frac{b}{\tau m}x_2-\left(\frac{1}{\tau}+\frac{b}{m}\right)x_3$, $g_0(x)=\frac{1}{\tau m}$, and $\omega=0.2\sin(t)$.  The cost functional for the  nominal system is selected as $r(x,u_0)=\|x\|^2+u_0^2$, and consequently the optimal value function and optimal control policy are given by $V^*(x)=2.2669x_1^2+2.3580x_2^2+0.0470x_3^2+3.1390x_1x_2+0.2x_1x_3+0.2534x_2x_3$ and $u_0^*(x)=x_1+1.2669x_2+0.4695x_3$, respectively. The compact set of interest is set as $\chi=[-1 ~1]\times [-1 ~1]\times [-5 ~5]$, and the data points to extrapolate the BE are selected on a $5\times 5\times 5$ data grid covering the domain $\chi$. The ESO is designed with $L=[4 ~6 ~4 ~1]^{\rm{T}}$ and $\varepsilon=0.01$. The basis function is selected as $\phi(x)=[x_1^2 ~x_2^2 ~x_3^2 ~x_1x_2 ~x_1x_3 ~x_2x_3]^{\rm{T}}$. The adaptation gains of the RL adaptive laws are the same as those in Example 1.

\begin{figure}
  \centering
  \includegraphics[width=0.40\textwidth,bb=15 0 320 188, clip]{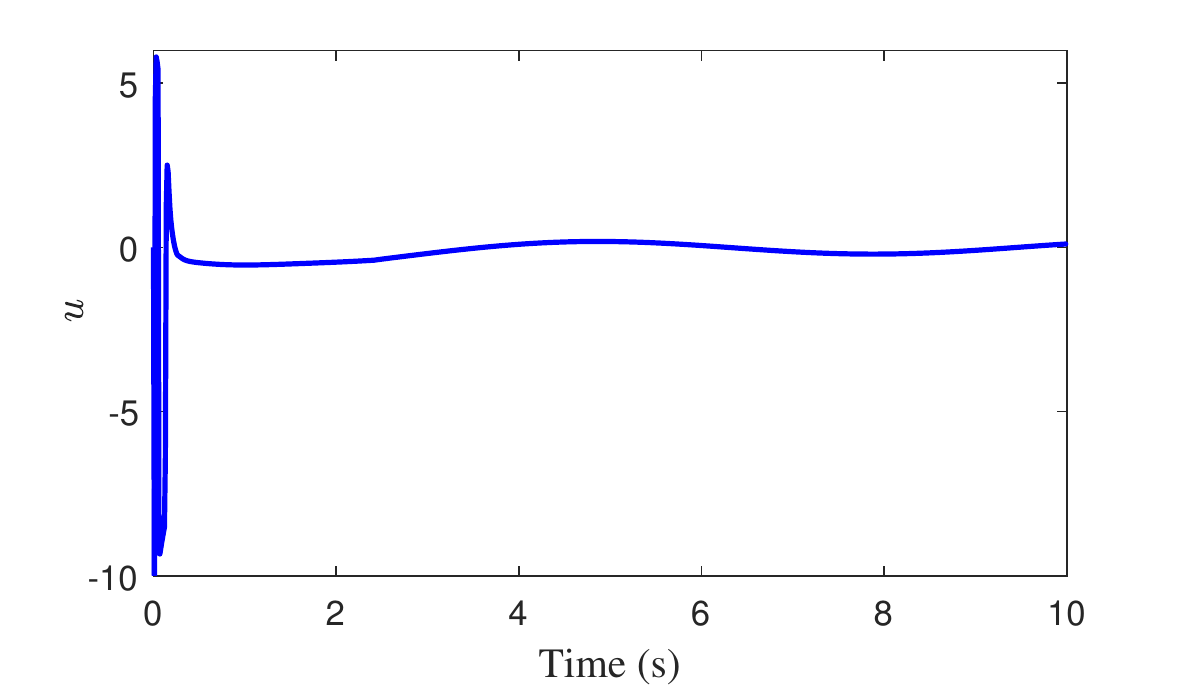}
  \caption{Trajectory of the control input $u$ (Example 2).}\label{fig13}
\end{figure}

\begin{figure}
  \centering
  \includegraphics[width=0.40\textwidth,bb=5 0 320 240, clip]{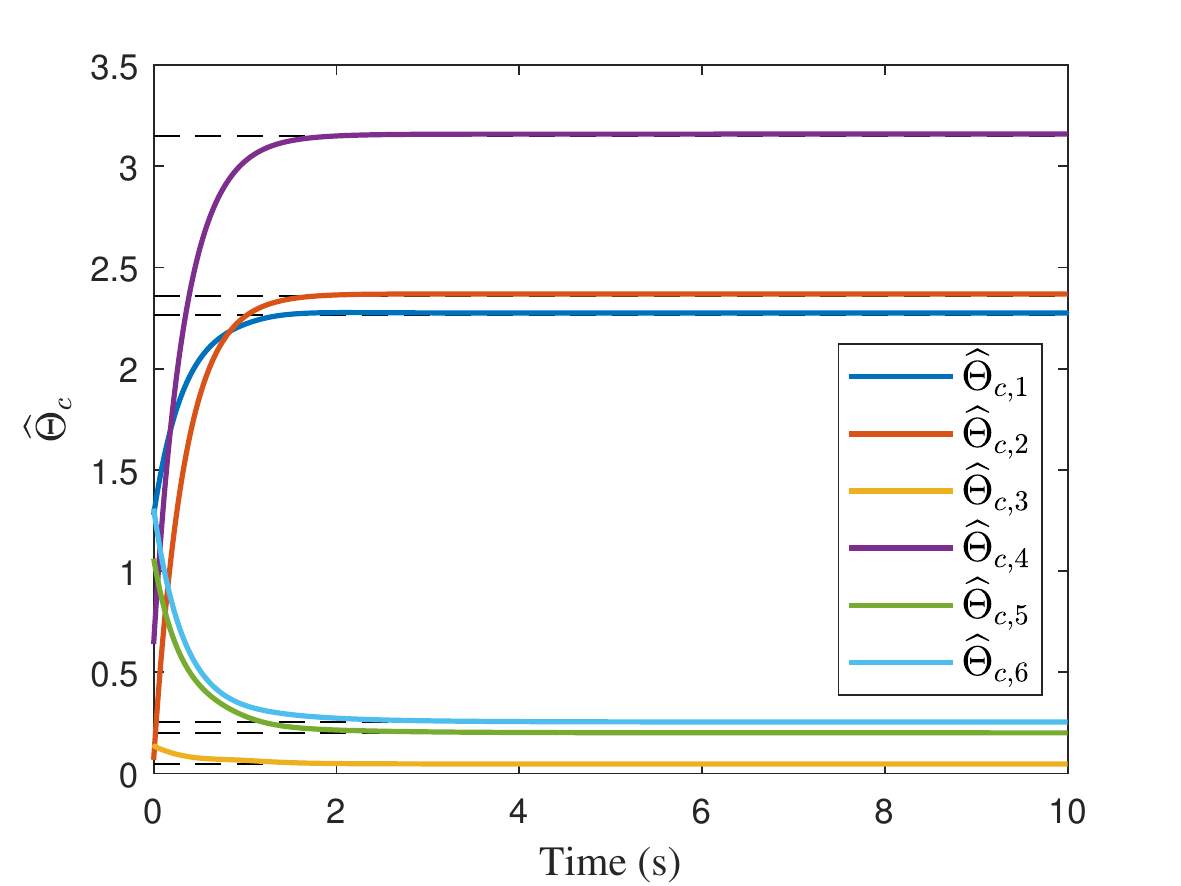}
  \caption{Trajectories of the actor NN weight (Example 2).}\label{fig14}
\end{figure}

Simulation is done with initial conditions $x(0)=[0.5 ~0.5 ~4]^{\rm{T}}$, $[\widehat{x}_1(0) ~\widehat{x}_2(0)~\widehat{x}_3(0)]^{\rm{T}}=[0 ~0 ~0]^{\rm{T}}$, $\Gamma(0)=\textrm{diag}\{100,100,100,100,100,100\}$, and the elements of $\widehat{\Theta}_v(0)$ and $\widehat{\Theta}_c(0)$ randomly selected on $[0 ~2]$. Figs. \ref{fig12}-\ref{fig14} show the simulation results, from which one can see that the developed approach achieves satisfactory performance for the third-order kinematic model.

\section{Conclusion}

An ESO-based reinforcement learning and disturbance rejection framework is established for uncertain nonlinear systems having non-simple nominal models.
The developed approach compensates for the total uncertainty and approximates the optimal policy for the compensated system simultaneously. Simulation of experience based RL is employed to utilize the nominal model and to relax the requirement of the PE condition. The obtained results provide a novel learning-based solution for the disturbance rejection of uncertain nonlinear systems, especially those having non-simple nominal models, which are quite common in practice. Future research works will be directed at the extension of the results to systems with constraints \cite{Yang-2020b,Yang-2020c}, and the application of the developed approach to robotic systems \cite{Li-2019}.

\end{document}